\documentclass[11pt,envcountsame]{llncs}

\newif\ifcomments
\commentstrue  

\usepackage[usenames,dvipsnames]{color}
\usepackage[utf8]{inputenc}
\usepackage[colorlinks=true]{hyperref}
\usepackage[square,sort,comma,numbers,sectionbib]{natbib}
\usepackage{amsmath, amssymb}
\usepackage{fullpage}
\usepackage{mathtools}
\usepackage[inline]{enumitem}
\usepackage{comment}
\usepackage{tikz,tikz-cd, forest}
\usetikzlibrary{hobby}
\usepackage{soul}
\usepackage[hyphenbreaks]{breakurl} 

\usepackage{xifthen}
\newcommand{\ifequals}[3]{\ifthenelse{\equal{#1}{#2}}{#3}{}}

\usepackage[ruled,vlined,linesnumbered]{algorithm2e}
\DontPrintSemicolon
\SetKwProg{Fn}{function}{}{}
\SetKwFunction{GCD}{GCD}
\SetKwFunction{Resultant}{Resultant}
\SetKwFunction{Roots}{Roots}

\newcommand{\C}{\mathbb{C}}
\newcommand{\F}{\mathbb{F}}
\newcommand{\FF}{\mathbb{F}}
\newcommand{\QQ}{\mathbb{Q}}
\newcommand{\ZZ}{\mathbb{Z}}

\newcommand{\Z}{\mathbb{Z}}

\newcommand{\PP}{\mathbb{P}}
\newcommand{\X}{\mathcal{X}}
\DeclareMathOperator{\End}{End}

\DeclareMathOperator{\Jac}{Jac}

\newcommand{\Ksing}{K^{\textrm{sing}}}
\newcommand{\Ksmooth}{K^{\textrm{sm}}}

\let\braket\relax
\DeclarePairedDelimiterX{\braket}[2]{\langle}{\rangle}{#1 \delimsize\vert #2}

\DeclareMathOperator{\poly}{poly}

\DeclarePairedDelimiter{\abs}{\lvert}{\rvert}
\let\ket\relax
\DeclarePairedDelimiter{\ket}{\lvert}{\rangle}

\def\C{\mathbb{C}}

\def\Z{\mathbb{Z}}
\def\F{\mathbb{F}}

\def\X{\mathcal{X}}

\spnewtheorem{heuristic}[proposition]{Heuristic}{\bfseries}{\itshape}

\definecolor{myyellow}{rgb}{1.0, 0.75, 0.0}
\definecolor{mygreen}{rgb}{0.35, 0.71, 0.0}

\ifcomments
\newcommand{\KS}[1]{\textcolor{Bittersweet}{{\sf (Kate:} {\sl{#1})}}}
\newcommand{\BO}[1]{\textcolor{red}{{\sf (Bo:} {\sl{#1})}}}
\newcommand{\BF}[1]{\textcolor{blue}{{\sf (Boris:} {\sl{#1})}}}
\newcommand{\JD}[1]{\textcolor{red}{{\sf (Jake:} {\sl{#1})}}}

\newcommand{\RB}[1]{\textcolor{olive}{{\sf (Ross:} {\sl{#1})}}}

\newcommand{\BS}[1]{\textcolor{magenta}{{\sf (Ben:} {\sl{#1})}}}

\newcommand{\xl}[1]{\textcolor{orange}{{\sf (Christelle:} {\sl{#1})}}}
\newcommand{\LZ}[1]{\textcolor{violet}{{\sf (Lukas:} {\sl{#1})}}}
\else
\newcommand{\KS}[1]{}
\newcommand{\BO}[1]{}
\newcommand{\BF}[1]{}
\newcommand{\JD}[1]{}
\newcommand{\RB}[1]{}
\newcommand{\BS}[1]{}
\newcommand{\xl}[1]{}
\newcommand{\LZ}[1]{}
\newcommend{\credit}[1]{\textcolor{yellow}{{\sf (Credit:} {\sl{#1})}}} 
\fi

\title{Failing to hash into supersingular isogeny graphs}

\author{
    Jeremy Booher\inst{1}
    \and
    Ross Bowden \inst{2}
    \and
    Javad Doliskani\inst{3}
    \and
    Tako Boris Fouotsa \inst{4}
    \and
    Steven D. Galbraith\inst{5}
    \and
    Sabrina Kunzweiler \inst{6}
    \and
    Simon-Philipp Merz\inst{7}
    \and
    Christophe Petit\inst{8,13}
    \and
    Benjamin Smith\inst{9}
    \and
    Katherine E. Stange\inst{10}
    \and
    Yan Bo Ti\inst{11}
    \and
    Christelle Vincent\inst{12}
    \and
    Jos\'e Felipe Voloch\inst{13}
    \and
    Charlotte Weitk{\"a}mper\inst{14}
    \and
    Lukas Zobernig\inst{5}
}

\institute{
    {Department of Mathematics,
    University of Florida, Gainsville, Florida, USA
    \\
    \email{jeremybooher@ufl.edu}}
    \and
    {Department of Computer Science, University of Bristol, Bristol, UK\\
    \email{ross.bowden@bristol.ac.uk}} 
    \and
    {Department of Computer and Software, 
    McMaster University, Hamilton, Canada 
    \\
    \email{jake.doliskani@mcmaster.ca}}
    \and
    {LASEC, EPFL, Lausanne, Switzerland \\
    \email{tako.fouotsa@epfl.ch}} 
    \and
    {Department of Mathematics,
    The University of Auckland, Auckland, New Zealand
    \\
\email{s.galbraith@auckland.ac.nz,lukas.zobernig@gmail.com}}
    \and 
    {Ruhr-Universit\"at Bochum, Bochum, Germany\\
    \email{sabrina.kunzweiler@ruhr-uni-bochum.de}}
    \and
    {Department of Computer Science, ETH Zurich, Switzerland \\
    \email{research@simon-philipp.com}}
    \and
    Laboratoire d'Informatique, \\
    Universit\'e libre de Bruxelles, Bruxelles, Belgium\\
    \email{christophe.f.petit@gmail.com}
    \and
    {Inria and Laboratoire d'Informatique (LIX),
    CNRS, École polytechnique, Institut Polytechnique de Paris,
    Palaiseau, France
    \\
    \email{smith@lix.polytechnique.fr}} 
    \and
    {Department of Mathematics,
    University of Colorado Boulder, Boulder, Colorado, USA}
        \\
        \email{kstange@math.colorado.edu}
    \and
    {DSO, Singapore\\ 
    \email{yanbo.ti@gmail.com}}
    \and
    {Department of Mathematics and Statistics, University of Vermont, Burlington, Vermont, USA\\
    \email{christelle.vincent@uvm.edu} }
    \and
    {School of Mathematics and Statistics,
    University of Canterbury, Christchurch, New Zealand
    \\
    \email{felipe.voloch@canterbury.ac.nz}}
    \and
    {University of Birmingham, Birmingham, UK \\
    \email{c.weitkaemper@pgr.bham.ac.uk}}
}

\begin{document}

\maketitle

\begin{abstract}
    An important open problem in supersingular isogeny-based cryptography is to produce, without a trusted authority, concrete examples of ``hard supersingular curves'' that is, equations for supersingular curves for which computing the endomorphism ring is as difficult as it is for random supersingular curves.  A related open problem is to produce a hash function to the vertices of the supersingular $\ell$-isogeny graph which does not reveal the endomorphism ring, or a path to a curve of known endomorphism ring.  Such a hash function would open up interesting cryptographic applications.  In this paper, we document a number of (thus far) failed attempts to solve this problem, in the hope that we may spur further research, and shed light on the challenges and obstacles to this endeavour.  The mathematical approaches contained in this article include: \begin{enumerate*}[label=(\roman*)] \item iterative root-finding for the supersingular polynomial; \item gcd's of specialized modular polynomials; \item using division polynomials to create small systems of equations; \item taking random walks in the isogeny graph of abelian surfaces, and applying Kummer surfaces; and \item using quantum random walks. \end{enumerate*} \textbf{Keywords:} isogeny-based cryptography, hashing, elliptic curves.
\end{abstract}

\begingroup
\makeatletter
\def\@thefnmark{$\ast$}\relax
\@footnotetext{\relax
\def\ymdtoday{\leavevmode\hbox{\the\year-\twodigits\month-\twodigits\day}}\def\twodigits#1{\ifnum#1<10 0\fi\the#1}%
Date of this document: \ymdtoday.
}
\endgroup

\section{Introduction}

Supersingular curves (and isogenies between them) have become a hot topic in cryptography over the last ten years or so. Fortunately the theory of complex multiplication provides efficient algorithms to generate a supersingular elliptic curve over $\FF_{p^2}$, even for the astronomically large $p$ that are used for cryptographic applications (see Br{\"o}ker~\cite{Bro2009}). It is also known how to uniformly sample a supersingular elliptic curve over $\FF_{p^2}$: generate one curve $E_0$ using Br{\"o}ker's method and then take a sufficiently long random walk in the supersingular isogeny graph to get a curve $E$.

There are several flavors of isogeny-based cryptography.
One of the earliest proposals was the cryptographic hash function based on isogenies, by Charles, Goren and Lauter~\cite{CGL}.
Another early proposal was to obtain a group action of the ideal class group on a set of elliptic curves. This was first proposed by Couveignes~\cite{Couv06} and re-discovered by Rostovtsev and Stolbunov~\cite{RS06}. Class group actions were made practical with the CSIDH scheme by Castryck, Lange, Martindale, Panny and Renes~\cite{CLMPR18}.
Various digital signature schemes have been proposed~\cite{YAJJS17,GPS20,DFG19,DPV19,BKV19,DFKLPW20}.
But the most studied isogeny-based cryptosystem of all is the key exchange protocol SIDH, by Jao and De Feo~\cite{jao2011towards}. Public keys in the SIDH protocol include not only elliptic curves but also certain auxiliary points on these curves.
The SIDH protocol has been a highly active area of research for over 10 years, but very recently major advances in cryptanalysis by Castryck and Decru~\cite{CD22}, Maino and Martindale~\cite{MM22} and Robert~\cite{Rob22} completely break SIDH, by exploiting the auxiliary points.
It remains to be seen whether some variant of SIDH can be secure and practical. Note that the other areas of isogeny-based cryptography, such as CSIDH and the signature schemes, do not use auxiliary points and so are not affected by the attack.
This paper was written before SIDH was broken, and we will refer to some results and papers that may not be relevant anymore.
Nevertheless, the general problems considered in this paper are still relevant for isogeny cryptography and remain worthy of study.

One of the main computational problems in isogeny-based cryptography  is to compute an isogeny between two given supersingular elliptic curves over the same finite field $\FF_{p^2}$. This problem is called the \emph{supersingular isogeny problem} or the \emph{path finding problem in  the supersingular isogeny graph}. It is believed to be hard, even for quantum computers.
A related problem is the \emph{supersingular endomorphism ring problem}: Given a supersingular elliptic curve $E$ over $\F_{p^2}$, compute its endomorphism ring $\End(E)$ (or even just one non-trivial endomorphism of $E$).
The supersingular endomorphism ring problem and the supersingular isogeny problem are related \cite{EisentragerHLMP18,GalbraithPST16,Wesolowski21}.

The algorithm using complex multiplication sketched in the first paragraph for generating a uniformly distributed supersingular curve has the side-effect that the person who generated the curve also knows a path from $E_0$ to $E$.
In certain cryptographic applications this approach is not acceptable as it allows a user to insert a trapdoor or in some other way violate the desired security.
There are a number of papers that have already mentioned this problem \cite{LB2020,CPV20,ADFMP20,de2019verifiable,MTTY20}. Currently the only solution known is to involve some ``trusted party'' to generate a random curve and then ``forget'' any resulting secret information.  See \cite{BassoEtAl} for trusted-setup solutions.
There is great interest in finding better ways to solve this problem that do not require trusting a single party. 
Among other applications, 
it would circumvent the trusted setup in an isogeny-based verifiable delay function~\cite{de2019verifiable}, in delay encryption~\cite{burdges2021delay} and in an SIDH-based oblivious pseudorandom function~\cite{boneh2020oblivious}. For the latter, the necessity of the trusted setup was pointed out by~\cite{basso2021cryptanalysis}.
Before SIDH was broken, using a starting curve that is generated uniformly at random would prevent torsion point attacks~\cite{DBLP:conf/asiacrypt/Petit17,quehen2021improved,kutas2021one}.

Applications of hashing to hard supersingular curves might include hash-and-sign signatures, oblivious pseudorandom functions~\cite{boneh2020oblivious} and password-authenticated key exchange~\cite{azarderakhsh2020not}.

There are (at least) three general problems that are of interest for isogeny-based cryptography:
\begin{enumerate}
    \item Given a prime $p$, to compute a supersingular curve $E$ over $\FF_{p^2}$ without revealing anything about the endomorphism ring or providing any information to help solve the isogeny problem (for isogenies from $E$ to some other supersingular curve over $\FF_{p^2}$).
    This is the problem of \textbf{demonstrating a hard curve}~\cite{LB2020}.
    
    \item\label{item:random-curve} Given a prime $p$, to generate \textbf{uniformly random} supersingular curves $E$ over $\FF_{p^2}$ without revealing anything about the endomorphism ring or providing any information to help solve the isogeny problem to other supersingular curves over $\FF_{p^2}$.
    
    \item \textbf{Defining a hash function to the entire supersingular graph}. To produce a hash function taking arbitrary strings as input, and outputting supersingular $j$-invariants. The hard problems in this context include both pre-image finding and collision-finding for the hash function, and path finding and endomorphism ring computation for the output curve.  We ask for these problems to remain hard on curves produced by the hash function.
\end{enumerate}

There are also variants of these problems that involve sampling from (resp. mapping to) subsets of the set of supersingular curves. The most significant is \textbf{defining a hash function just to the $\FF_p$ subgraph}.

The two obvious approaches to these problems are to use tools from the theory of complex multiplication and/or random walks. However neither method is secure for our problems. The insecurity of methods based on random walks is self-evident. The insecurity of methods based on CM is less clear, and was demonstrated by 
Castryck, Panny and Vercauteren~\cite{CPV20} and Love and Boneh~\cite{LB2020}.
We refer to Section~\ref{sec:sec2} for details.

Castryck--Panny--Vercauteren~\cite{CPV20} and Wesolowski~\cite{Wesolowski21Orientations} have considered the analogous approach in the special case of sampling supersingular curves with $j$-invariant in $\FF_p$ using CM theory.
Again they show that any such approach is not secure (they show how to solve the class group action problem in subexponential and polynomial time respectively). 

Hence we need new ideas. The goal of the paper is to explain some possible approaches and to discuss the obstructions to getting a practical solution.

In all cases we are interested in an efficient  algorithm that takes as input $p$, can be executed without any secret information, and that outputs (the $j$-invariant of) a supersingular elliptic curve over $\FF_{p^2}$. We do not want the algorithm to provide any additional information that would be useful to the person who executes it.
For the problem of generating a single hard curve (e.g., to bypass the requirement for trusted set up), the meaning of ``efficient'' might be relaxed, as long as it is feasible in applications.

As already mentioned, it would already be interesting to have an algorithm that returns a single curve. But the most desirable outcome is a \textbf{cryptographic hash function} $H(m)$ that takes a binary string $m$ and returns a supersingular $j$-invariant and satisfies these properties:
    \begin{enumerate}
    \item It is efficient and deterministic.
    \item It is hard to find a \textbf{collision}, namely two binary strings $m_1$ and $m_2$ such that $H(m_1)  = H(m_2)$.
    \item It is hard to \textbf{invert}, namely given an $h$ in the codomain it is hard to compute a binary string $m$ such that $H(m) = h$.
    \item The $j$-invariants are uniformly distributed in the codomain.
    \end{enumerate}

Note that one can build an algorithm for hashing to the supersingular set by combining a standard cryptographic hash function $H'$ (e.g., SHA-3) with a randomised algorithm to generate a supersingular curve (as in problem 2 listed above). To do this, simply compute $H'(m)$ and use it as the seed to a pseudorandom generator and then run the algorithm to generate a supersingular curve replacing all calls to randomness with this pseudorandom sequence.
Hence, it suffices to focus on problems 1 and 2 above.


Several of the approaches in our paper try to bypass the problem of working with polynomials of exponentially-large degree.
Section~\ref{sec:Iterating} sketches an approach motivated by iterated methods for root-finding (such as the Newton-Raphson method). However the main idea in this section is to avoid writing down the polynomial by indirectly computing its evaluation at a given point. This motivates a study of iterative methods in this special case.
Similarly, Section~\ref{sec:modular} studies an approach based on modular curves and the fact that one can compute the roots in $\F_{p^2}$ of the greatest common divisor of two polynomials $F(x,x^p)$ and $G(x,x^p)$ in polynomial time in certain circumstances, even though the polynomials themselves have exponential degree.
This approach does not lead to a useful solution at present, as the computation only produces curves that could feasibly have been computed using the CM method.
Section~\ref{sec:reverse} also attempts to control the growth of polynomials, by giving a system of low-degree polynomials whose common solution would give a desired curve.

Other methods try to use random walks in new ways.
Section~\ref{sec:genus_2} suggests walking on the isogeny graph of abelian surfaces, until one lands on a reducible surface.
The challenge faced by this method is that reducible surfaces are exponentially rare in the isogeny graph and we lack techniques to navigate to one from an arbitrary position in the graph.
Finally, Section~\ref{sec:quantum} suggests a way to use a quantum analog of the CGL hash to generate a random supersingular curve.  The way a quantum algorithm uses randomness means this cannot be combined with a standard cryptographic hash function as described above.  If properly implemented on a quantum computer, the algorithm makes the path information inaccessible to the user. But without a method to certify the use of the quantum algorithm, this approach only replaces the need for a trusted entity from one who will erase the path data to one who will promise to use a quantum computer.

Between release and revisions for this work, the concurrent work \cite{hashingalso}, which also proposes some approaches to the hashing problem, was made public.  In particular, the papers appear to overlap in suggesting a system of equations based on torsion point restrictions (compare Section~\ref{sec:rev-system} and \cite[Section 6.3]{hashingalso}).

We hope the ideas and analysis in our paper will be useful to researchers. We identify a number of obstructions to efficient hashing to supersingular curves. We hope that future research might overcome one of these obstructions.

\subsection*{Funding}
This work was supported by the 
Marsden Fund Council administered by the Royal Society of New Zealand [to J.B.];
Natural Sciences and Engineering Research Council of Canada (NSERC) [to J.D.];
Deutsche Forschungsgemeinschaft (DFG, German Research Foundation) under Germany's Excellence Strategy [EXC 2092 CASA - 390781972 to S.K.];
Engineering \& Physical Sciences Research Council [EP/P009301/1 to S.-P.M., EP/S01361X/1 and EP/V011324/1 to C.P., EP/T517872/1 to R.B.];
l'Agence nationale de la recherche (ANR) [program CIAO (ANR-19-CE48-0008) and a \emph{Plan France 2030} grant (ANR-22-PETQ-0008 PQ-TLS) to B.S.];
National Science Foundation [NSF-CAREER CNS-1652238 to K.S., DMS-1802323 to C.V.];
Simons Foundation [Fellowship 822143 to K.S.];
Ministry of Business, Innovation and Employment and the Marsden Fund Council administered by the Royal Society of New Zealand [to J.F.V., to S.G.].

\subsection*{Data availability statement}
\hl{No data available.}

\subsection*{Acknowledgements}
This project was initiated as part of the Banff International Research Station (BIRS) Workshop 21w5229, \emph{Supersingular Isogeny Graphs in Cryptography}.  The project owes a debt of gratitude to BIRS and to the organizers of that workshop: Victoria de Quehen, Kristin Lauter, Chloe Martindale, and Christophe Petit.  The project was led by Steven Galbraith, Christophe Petit, Yan Bo Ti, and Katherine E. Stange.
We would also like to thank Chloe Martindale for useful discussions, Annamaria Iezzi for her involvement in Section~\ref{sec:modular}, as well as Wouter Castryck and Eyal Goren for contributing ideas to Section~\ref{sec:genus_2}.

\section{Existing methods}
\label{sec:sec2}

We briefly review the existing methods of generating supersingular curves.  Neither is secure in the sense of the introduction.  Nevertheless, these two paradigms form the basis of the methods proposed in this work, which fall broadly into methods based on random walks, and methods based on finding roots to high degree polynomials (or systems of such).

The Charles-Goren-Lauter hash function \cite{CGL} hashes into the supersingular curves over $\FF_{p^2}$.  At each vertex of the supersingular isogeny graph, the out-directed edges are labelled in some fixed deterministic manner.  Starting from a known curve such as $j=1728$, the bitstring to be hashed is interpreted as directions for a walk through the graph, via the labelling just mentioned.  If the walk is sufficiently long, it is known from the properties of the graph (it is Ramanujan) that the endpoint will be uniformly randomly chosen from amongst all the vertices of the graph~\cite{CGL}.  However, the walk itself is a path to $j=1728$ and therefore the path-finding problem from the endpoint is trivial, unless this information is discarded by a trusted authority.

The CM method of Br\"oker \cite{Bro2009} finds supersingular roots of a Hilbert class polynomial.  The Hilbert class polynomial $H_{\mathcal{O},p}$ for a quadratic order $\mathcal{O}$ modulo $p$ is a polynomial in $\FF_p[x]$ whose roots in $\overline{\mathbb{F}}_p$ are the $j$-invariants of elliptic curves whose endomorphism rings contain a copy of $\mathcal{O}$.  In order to apply known root-finding algorithms, or indeed, to obtain the polynomial at all, the degree of $H_{\mathcal{O},p}$ must be small. 
This implies that $\mathcal{O}$ itself has non-integral elements of small norm.  The images of such elements in the endomorphism ring are termed \emph{small endomorphisms}, and so all the curves obtained have small endomorphisms. 
At the very least, it follows that the curves obtained are far from uniformly distributed.
Furthermore, having a small endomorphism is known to be a serious vulnerability \cite{LB2020,CPV20}. 
Precisely, Castryck, Panny and Vercauteren~\cite{CPV20} study the CSIDH case and show how to efficiently compute an ideal class as required to break CSIDH when given a small degree endomorphism.
Love and Boneh~\cite{LB2020} consider the more general case, such as arises in SIDH, and also show a general and efficient approach to computing isogenies between any two such curves in this setting.  In general, anything which reveals the endomorphism ring will be a vulnerability \cite{Wesolowski21}.

\section{Iterating to supersingular $j$-invariants} \label{sec:Iterating}

In this section, we propose a method for generating a hard curve.  
For a prime number $p>2$, define the polynomial $H_p(t)$, known as the \emph{Hasse polynomial} or \emph{supersingular polynomial}, by
\begin{equation} \label{eqn:Supersingular poly}
    H_p(t) = \sum_{j=0}^{(p-1)/2} \begin{pmatrix}\frac{p-1}{2}\\ j\end{pmatrix}^2 t^j.
\end{equation}

\begin{proposition}
Let $E_\lambda$ denote the elliptic curve whose Legendre form is $y^2 = x(x-1)(x-\lambda)$.  Then for $\lambda \in \F_p$
\[
\# E_\lambda(\F_p) \equiv p + 1 - H_p(\lambda) \pmod{p}.
\]
Similarly for $\lambda \in \F_{p^2}$
\[
\# E_\lambda(\F_{p^2}) \equiv p^2 + 1 - H_p(\lambda)^{p+1} \pmod{p}.
\]
\end{proposition}

\begin{proof}
This follows from the proof of \cite[Theorem V.4.1(b)]{Silverman}. 
\qed
\end{proof}

Thus $\lambda\in \F_{p^2}$ is a root of $H_p(t)$ if and only if $E_\lambda$ is a supersingular elliptic curve.   It is known that all the roots belong to $\F_{p^2}$ and that, for $p \equiv 3 \pmod 4$, we have that $p^{1/2+o(1)}$
of them belong to $\F_p$.  None belong to $\F_p$ when $p \equiv 1 \pmod{4}$. 
This follows since the number of supersingular curves over $\F_p$ is 
$p^{1/2+o(1)}$ by combining
\cite{Goldfeld} and \cite[Eq (1)]{Delfs-Galbraith}
and such a curve
can be put in Legendre form if and only if all of its $2$-torsion is rational, which is only possible when 
$p \equiv 3 \pmod 4$.

The basic idea is to compute a random root of the polynomial, thus giving a random supersingular elliptic curve.  At first glance this seems impractical, as representing the polynomial $H_p(t)$ takes exponential space, and computing $H_p(t)$ would take exponential time.
However, we can compute $H_p( \lambda )$ for $\lambda \in \FF_p$ in polynomial time using Schoof's algorithm to compute $\# E_\lambda(\F_p)$.  For $\lambda \in \F_{p^2}$, we can similarly compute $H_p(\lambda)^{p+1}$ by computing $\# E_\lambda(\F_{p^2})$.
It is unclear whether there is a fast way to compute $H_p( \lambda )$ for $ \lambda \in \F_{p^2}$. 

\subsection{Iterating to a root} \label{ss:iteratingtoroot} One approach to finding a root of $H_p(t)$ is to iterate a polynomial function over a finite field as inspired by the Newton-Raphson method.
Recall that the Newton-Raphson method finds a root of a polynomial $f(x)$ by first picking a point on the domain $t_0$ and iteratively computing
\begin{equation} \label{iteration0}
    t_{n+1} = t_n - \frac{f(t_n)}{f'(t_n)}
\end{equation}
(while $f'(t_n) \neq 0$).  
If a fixed point $t_{m+1}=t_m$ is found, then we can conclude that $f(t_m)=0$ and that we have found a root.

In this vein, our ``preliminary'' idea is to find the roots using the same method. 
So one picks some $t_0\in\F_p$ (or $\F_{p^2}$), and then defines
\begin{equation}\label{iteration1}
t_{n+1}=t_n - \frac{H_p(t_n)}{H_p'(t_n)}
\end{equation}
(assuming $H_p'(t_n) \neq 0$).
It is clear that if $t_{m+1}=t_m$, we must have that $H_p(t_m)=0$, and we have found a supersingular elliptic curve.
Furthermore, this method could allow us to define a hash function into supersingular curves, by using the hash input to determine $t_0$ and then iterating \eqref{iteration1}. 

However, there are three issues with this idea:
\begin{enumerate}
    \item The algorithm may not halt at a fixed point (the iteration may become stuck in a cycle).
    \item The algorithm may reach a fixed point, but require too many iterations to efficiently compute.
    \item We do not know how to compute $H_p'(t)$ efficiently, or compute $H_p(t)$ efficiently for $t \in \F_{p^2} \setminus \F_p$.  
\end{enumerate}

To eliminate the third obstacle, we can consider the following alternatives to  the Newton-Raphson method which share the key property that fixed points of the iteration correspond to roots of $H_p(t)$:
\begin{align}
    t_{n+1} &=t_n - H_p(t_n) \label{iteration2} \\
    t_{n+1} &=t_n - H_p(t_n)^{p+1}\, . \label{iteration3}
\end{align}
The denominator $H_p'(t)$ in the previous attempt speeds up convergence to a root in a field with a metric, if we are already close to a root. 
In a finite field with the discrete topology there is no reason to include the denominator,  and removing the denominator also removes the possibility the iteration is undefined.
Using Schoof's algorithm we may efficiently iterate \eqref{iteration2} over $\F_{p}$ (which is only of interest when $H_p(t)$ has roots over $\F_p$, i.e. $p \equiv 3 \pmod{4}$), while we may efficiently iterate \eqref{iteration3} over $\F_{p^2}$.

The first two obstacles are thornier to tackle: there are plenty of choices of $t_0$ where iteration leads to a cycle and not a fixed point, and paths to a fixed point can be very long.  
These are fundamental obstructions which we will discuss experimentally and compare to the behavior of random mappings.


\subsection{Does Iteration Mimic Iteration of a Random Function}

\label{ss:iterationrandom}
In terms of understanding whether these iterative methods are useful for finding supersingular elliptic curves, the important quantities to understand for the iteration are:
\begin{enumerate}
    \item  the number of fixed points;
    \item  the number of points which eventually reach a fixed point upon iteration;
    \item  for these points, the maximum number of iterations needed to reach the fixed point; and
    \item  the number of points which reach a fixed point after $k$ iterations.
\end{enumerate}
We will be mainly interested in understanding to what extent these iterative methods look like iterating a random function (which we can understand theoretically).

Consider a random function from a set $S$ of size $n$ to itself.  In other words, the image of each element of $S$ is chosen independently and uniformly at random from $S$.  
The expected number of fixed points for a random function is one, and given our knowledge about iterating random functions we would not expect iteration to quickly close in on the fixed point.
However, the function we are iterating has many fixed points, and it is not a priori clear how iteration will behave.  

\begin{heuristic} \label{heuristic:random}
Iterating \eqref{iteration1} (resp. \eqref{iteration2}, \eqref{iteration3}) over $\F_q$ when $n=p, p^2$ (resp. when $n=p \equiv 3 \mod{4}$, when $n = p^2$) behaves like iterating a random function from a set of size $n$ to itself that has $\sqrt{n}$ fixed points.  
\end{heuristic}

We will experimentally explore the behavior of the iterative methods in this subsection, and theoretically study the behavior of iterating random functions with many fixed points in Section~\ref{ss:manyfixed}.

Experimentally, it appears that ``many'' points eventually reach a fixed point after iteration, which means that our iterative methods have a reasonable chance of finding a supersingular curve.
However, the maximum number of iterations needed seems to be on the order of $\sqrt{n}$, which is too long to be practical.  This is in line with the expected ``tail length'' of a random mapping \cite[Theorem 8.4.8]{mpq19}.  (Section 8.4 of \emph{loc. cit.} contains a survey of the properties of random mappings.)  Finally, the number of points which reach one of the $m$ fixed points after $k$ iterations appears to be on the order of  $(k+1) m$, at least when $k$ is small relative to $n$.  This matches the conclusions of our analysis of random functions with many fixed points in Section~\ref{ss:manyfixed}.   

While we provide examples of all the iterative methods, we have focused on iteration \eqref{iteration2} over $\F_p$ when $p \equiv 3 \pmod{4}$ as it is efficiently computable and well-motivated by analogy with the Newton-Raphson method.  The behavior we are seeing does not seem sensitive to the exact iterative method used. 

\begin{example}
When using the original Newton iteration \eqref{iteration1}, many points eventually reach a fixed point upon iteration.  For example, when $p=101$ the polynomial $H_p(t)$ has $50$ roots all defined over $\F_{101^2}$.  If one iterates using \eqref{iteration1}, $328$ of the elements of $\F_{101^2}$ eventually end up at a fixed point (about three percent).  In those cases it took at most $10$ iterations to reach a fixed point.  Similarly, when $p=211$ around twenty eight percent ($12747$ out of $211^2$ of the elements) eventually reach a fixed point.  The maximum number of iterations needed was $90$.  When $p=1009$, $800$ out of $10,000$ randomly chosen elements of $\F_{p^2}$ eventually reached a fixed point. 
\end{example}


\begin{example}
The behavior when iterating using \eqref{iteration2} over $\F_{p^2}$ is broadly similar; removing division by $H_p'(t_n)$ does not seem to have a significant effect.  For example when $p = 1009$,  $278$ of $1000$ randomly chosen elements of $\F_{1009^2}$ ended in a fixed point. More systematically, if we look at all primes between $30$ and $200$ and compute the percentage of elements of $\F_{p^2}$ which eventually reach a fixed point, the minimum and maximum percentages are about $1.9$ percent and $81$ percent.  The mean is about $30$ percent.  There are often quite long paths which eventually lead to a fixed point.
\end{example}

\begin{example}
Iterating using \eqref{iteration2} over $\F_p$ is only interesting when there are fixed points defined over $\F_p$, i.e. when $p \equiv 3 \pmod{4}$.  It appears broadly similar to the previous iterations considered.  
The number of fixed points is $p^{1/2+o(1)}$. 
Experimentally it looks like a sizeable fraction of the points of $\F_p$ eventually reach a fixed point, and that for small $k$ the number of points which reach a fixed point after $k$ iterations is about $(k+1)$ times the number of fixed points (so on the order of $(k+1) \sqrt{p}$).  The largest number of iterations needed to reach a fixed point appears to be on the order of $\sqrt{p}$.

To quantify this, we computed the minimum and maximum values for:
\begin{itemize}
    \item the number of fixed points divided by $\sqrt{p}$, denoted $F_1(p)$;
    \item   the number of elements of $\F_p$ iterating to a fixed point, divided by $p$, denoted $F_2(p)$;
    \item   the largest number of iterations needed to reach a fixed point divided by $\sqrt{p}$, denoted $F_3(p)$.
\end{itemize}
Table~\ref{table:iteratingdata} shows the minimum and maximum values of these values for primes in several ranges.

\begin{table}[ht]
\centering
\begin{tabular}{c|c c| c c| c c }
$p$ in Range: & $\min F_1(p)$ & $\max F_1(p)$ & $\min F_2(p)$ & $\max F_2(p)$ & $\min F_3(p)$ & $\max F_3(p)$ \\ 
\hline
$100$ to $2000$ & $.23$ & $3.9$ & $.019$ & $.93$ & $.034$ & $.95$ \\
$2000$ to $3000$ & $ .29$ & $4.0$ & $.014$ & $.61$ & $.062$ & $.64$\\
$20000$ to $21000$ & $.27$ & $4.0$ & $.0085$ & $.46$ & $.035$  & $.47$\\
\hline
\end{tabular}
\caption{Statistics about iteration \eqref{iteration2} over $\F_p$}
\label{table:iteratingdata}
\end{table}

Figure~\ref{fig:experiment-overp} shows a graph of the ratio of the number of elements of $\F_{p}$ which reach a fixed point after $5$ iterations of \eqref{iteration2} and of the number of fixed points, versus $p$.  As expected, this appears to be around $6$ but is somewhat noisy. 

\begin{figure}
    \centering
    \includegraphics[width = .75 \textwidth]{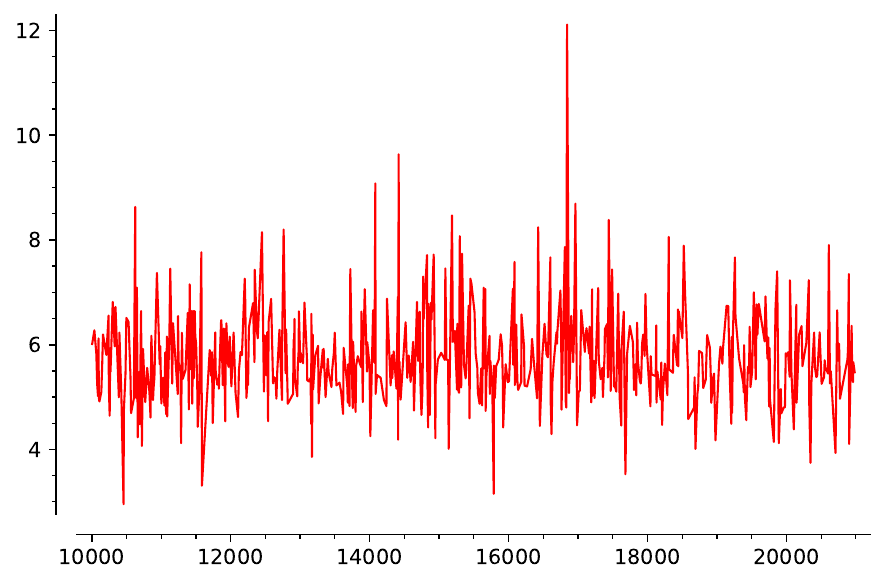}
    \caption{Number of elements which reach a fixed point after $5$ iterations of $x \mapsto x - H_p(x)$ divided by the number of fixed points, versus $p$}
    \label{fig:experiment-overp}
\end{figure}
\end{example}

\begin{example}
Iterating using \eqref{iteration3} over $\F_{p^2}$ preserves the cosets of $\F_p$ in $\F_{p^2}$.  At first glance, it looks like for most cosets the map behaves like a random map: each coset has very few fixed points and about $\sqrt{p}$ points in each coset lead to the fixed points.
\end{example}

\begin{example}
For a hundred randomly chosen function $\F_{101^2} \to \F_{101^2}$ with $50$ fixed points:
\begin{itemize}
    \item  the number of elements iterating to a fixed point ranged from $309$ to $10194$ ($3\%$ to $99\%$), with a mean of $4525$ ($44\%$).
    \item  the maximum number of iterations needed to reach a fixed point ranged from $12$ to $350$, with a mean of $116$.
    \item  The number of elements within $5$ iterations of a fixed point ranged from $181$ to $443$, with a mean of $295$.
\end{itemize}
This is broadly in line with the behavior seen in the previous examples, so experimentally it seems reasonable that iterating \eqref{iteration1} \eqref{iteration2} or \eqref{iteration3} behaves like iterating a random function, supporting Heuristic~\ref{heuristic:random}.
\end{example}

Based on this behavior, using iteration to efficiently find a supersingular elliptic curve (or to hash to a supersingular elliptic curve) does not seem to be practical.  For concreteness, we will focus on iterating \eqref{iteration2} over $\F_p$ when $p \equiv 3 \pmod{4}$.  The basic idea would be to pick a random starting element and iterate $k$ times, hoping to find a fixed point.  
Thus the key property of the iteration to understand is the number of points which iterate to a fixed point in $k$ steps.  We expect on the order of $(k+1) \sqrt{p}$ points with this property, both based on experiments and based on
 Heuristic~\ref{heuristic:random} plus the analysis in Section~\ref{ss:manyfixed}.

\begin{proposition}
Assuming Heuristic~\ref{heuristic:random}, iterating $k$ times over $\F_p$ will find a fixed point (and hence a supersingular elliptic curve) with probability on the order of $(k+1) / \sqrt{p}$.  This requires $k+1$ evaluations of $H_p(t)$.   
In particular, iteration would not give an efficient method of finding supersingular $j$-invariants.
\end{proposition}

\begin{proof}
   There are on the order of $(k+1) \sqrt{p}$ elements of $\F_p$ which iterate to a fixed point within $k$ steps.  To check whether the last element is fixed requires one additional evaluation of $H_p(t)$.  The chance of randomly choosing to start at one of these elements is on the order of $(k+1)/ \sqrt{p}$.   

While we can evaluate $H_p(t)$ efficiently using Schoof's algorithm, to be efficient the number of evaluations $k+1$ must still be polynomial in $\log(p)$. In that case the probability of finding a supersingular curve would be exponentially small in $\log(p)$. \qed
\end{proof}

\begin{remark}
   The probability of a random curve being supersingular is on the order of $1/\sqrt{p}$.  We can check whether a randomly chosen $j$-invariant is supersingular by evaluating $H_p(t)$.  Checking $k+1$ random curves would require $k+1$ evaluations of $H_p(t)$, and again would find a supersingular $j$-invariant in $\F_p$ with probability on the order of $(k+1)/\sqrt{p}$.  In particular, the iterative method is no better than randomly guessing and neither approach can efficiently find a supersingular $j$-invariant.
\end{remark}

\begin{remark}
For the iterative method to offer an improvement, we would need a way to make a ``giant step'' and efficiently iterate multiple times at once.  For example, given the $n$th iteration we would like to be able to efficiently  compute the $2n$th iteration.  We do not know if this is possible.
\end{remark}

\begin{remark}
This iterative method would not produce supersingular curves uniformly at random.
Taking $p=1019$, when iterating \eqref{iteration2} there are fixed points which no other elements of $\F_p$ reach upon iteration but there is also a fixed points that $50$ other elements of $\F_p$ reach upon iteration.
\end{remark}

\subsection{Random functions with fixed points} \label{ss:manyfixed}

We use the functional graph perspective on random mappings and the asymptotic analysis developed in Flajolet and Odlyzko~\cite{FO90} to analyze functions with many fixed points.   For a function $f: S \to S$, the vertices are $S$ and there is a directed edge from $a$ to $b$ if $f(a) = b$.  A function on $n$ elements with $m$ fixed points is represented by:
\begin{itemize}
\item  A functional graph, consisting of $m$ rooted trees (one for each fixed point) plus a set of components without fixed points;
\item  Each component without fixed points is a collection of (at least two) trees where the roots are permuted cyclically; 
\item  Each rooted tree is a node (the root) together with a possibly empty set of rooted trees that are the children.
\end{itemize}
Note that all of these objects are labeled.  

\begin{figure}
\[\begin{tikzcd}
	&& 13 & 16 &&& 3 &&& 12 \\
	11 & 8 && 9 &&& 14 && 4 \\
	&& 1 & 10 && 2 && 5 & 6 \\
	& 0 &&&&& 7 & 15
	\arrow[from=2-2, to=3-3]
	\arrow[from=4-2, to=3-3]
	\arrow[from=1-3, to=2-4]
	\arrow[from=3-3, to=3-4]
	\arrow[from=2-4, to=3-4]
	\arrow[from=2-1, to=2-2]
	\arrow[from=3-6, to=2-7]
	\arrow[from=2-7, to=3-8]
	\arrow[from=3-8, to=4-7]
	\arrow[from=4-7, to=3-6]
	\arrow[from=1-7, to=2-7]
	\arrow[from=2-9, to=3-8]
	\arrow[from=3-9, to=3-8]
	\arrow[from=1-10, to=2-9]
	\arrow[from=4-8, to=4-7]
	\arrow[from=1-4, to=2-4]
\end{tikzcd}\]
\caption{A functional graph on $17$ elements consisting of one rooted tree (with fixed point $10$) and one component based around a cycle of length $4$.}
\end{figure}

As in \cite{FO90}, there is a standard method to give relationships between exponential generating functions for these objects. 
Let $F_m(z)$ be the exponential generating function for random mappings with exactly $m$ fixed points.  This means that
\[
F_m(z) = \sum_{n=0}^\infty F_{m,n} z^n / n!
\]
where $F_{m,n}$ is the number of such functions with $n$ total elements.  Equivalently, it is the sum of $z^{|\varphi|}/|\varphi|!$ over all functions $\varphi$ with $m$ fixed points.
 Likewise let $C(z)$ and $T(z)$ be the exponential generating functions for components and trees, and let $C_{fpf}(z)$ be the exponential generating function for fixed-point-free components. 
 
 \begin{lemma}
 We have the following relationships:
\begin{align}
F_m(z) &= T(z)^m  \exp(C_{fpf}(z))\\
C_{fpf}(z) &= -\log(1-T(z)) - T(z) \\
T(z) &= z \exp(T(z)).
\end{align}    
 \end{lemma}
 
 \begin{proof}
 The first is a consequence of the fact that a function with $m$ fixed points consists of $m$ rooted trees plus a set of components with no fixed points.  
    It is standard that
 \[
 C(z) = \sum_{k \geq 1} \frac{1}{k} T(z)^k = \log(1/(1-T(z)))
 \]
 as a connected component based on a cycle of length $k$ is built out of $k$ trees and one can cyclically permute them.    Therefore we see that
 \[
 C_{fpf}(z) = C(z) - T(z) = \log(1/(1-T(z))) - T(z).
 \]
 The third is standard, a consequence of the fact that  a tree is a node plus a set of trees.  \qed
 \end{proof}
 
 We can use asymptotic analysis to compute the number of random functions  with $m$ fixed points.  Flajolet and Odlyzko
\cite[Proposition 1]{FO90} give an asymptotic expansion 
\[
T(z) = 1 - \sqrt{2} \sqrt{1-ez} - 1/3 (1-ez) + O\left((1-ez)^{3/2}\right)
\]
of $T(z)$ around its singularity at $z = 1/e$.  We can rewrite $F_m(z) = T(z)^m \exp(C_{fpf}(z))$ in terms of $T(z)$ as 
\[
F_m(z) = T(z)^m \frac{1}{1-T(z) } \exp(-T(z)) = T(z)^{m-1} \frac{z}{1-T(z)}
\]
using that $T(z) = z \exp(T(z))$.  
We have that 
 $z = \frac{1}{e} - \frac{1}{e}(1-ez)$, so the 
leading term in the asymptotic expansion of $F_m(z)$ is
\[
(e \sqrt{2}  \sqrt{1-ez})^{-1}.
\]
Using \cite[Theorem 1]{FO90} gives the asymptotic 
\begin{equation} \label{eq:asymptoticfmn}
\frac{F_{m,n}}{n!} \sim \frac{e^{n-1}}{\sqrt{2 \pi n}}.
\end{equation}
For example, taking $m=0$ gives the asymptotic $\frac{1}{e} \frac{1}{\sqrt{2 \pi n}} e^n$ for $F_{0,n}/n!$.  In comparison, Flajolet and Odlyzko's asymptotic analysis gave the known fact 
(letting $F_n$ denote the number of functions on $n$ elements) that $F_n / n! \sim \frac{1}{\sqrt{2 \pi n}} e^n $.  As expected, this implies that 
about $\frac{1}{e}$ of randomly chosen functions do not have a fixed point. 

\begin{remark}
Note that if $m$ is fixed as $n \to \infty$, the precise value of $m$ has no effect on the asymptotics of $F_{m,n}$.
\end{remark}

We will now modify our generating functions to take into account the number of elements which reach a fixed point after $k$ iterations of a random function.
The key case is for trees, where we consider the exponential generating function $T_k(z,u)$ where
the coefficient of $z^n u^\ell$ is the number of rooted trees of size $n$ with $\ell$ nodes that are distance at most $k$ from the root, divided by $n!$.

\begin{lemma}
 We have that $T_0(z,u) = u T(z)$,  that $T_k(z,u) = zu \exp(T_{k-1}(z,u))$, and that $T_k(z,1) = T(z)$.  
\end{lemma}

\begin{proof}
The first equality reflects that each tree has exactly one node at distance $0$ from the root.  The second comes from the fact that a rooted tree is a root plus a collection of child trees, and a node has distance at most $ k$ from the root if it either is the root or has distance at most $ k-1$ from the root of one of the child trees.  The third equality is clear.   \qed
\end{proof}

We likewise modify $F_m(z)$ to become a bi-variate exponential generating function $F_{m,k}(z,u)$ which counts nodes with distance at most $k$ to one of the $m$ fixed points.  It satisfies
\[
F_{m,k}(z,u) = T_k(z,u)^m \exp(C_{fpf}(z)) = T_k(z,u)^m \frac{z}{(1-T(z)) T(z)}.
\]

\begin{lemma}
The exponential generating function for the sum of the number of elements which reach a fixed point after $k$ iterations for functions with $m$ fixed points is 
\[
\frac{d F_{m,k}(z,u)}{du}|_{u=1}.
\]
\end{lemma}

\begin{proof}
This follows from viewing the bi-variate exponential generating function as a sum over all functions with $m$ fixed points. \qed
\end{proof}

\begin{proposition} \label{prop:expectedfixedpoints}
For fixed $m$ and $k$, the number  of elements which reach a fixed point after $k$ iterations for a random function on $n$ elements with $m$ fixed points is asymptotically $(k+1)m$ as $n \to \infty$.
\end{proposition}

\begin{proof}
We compute that $\frac{T_0(z,u)}{du}|_{u=1} = T(z)$ and
\begin{align*}
\frac{d T_{k}(z,u)}{du}|_{u=1} &= z \exp(T_{k-1}(z,1)) + z \exp(T_{k-1}(z,1)) \frac{d T_{k-1}(z,u)}{du}|_{u=1}\\
&= T(z) + T(z) \frac{d T_{k-1}(z,u)}{du}|_{u=1}\\
&= T(z)+ T(z)^2 + \ldots +T(z)^{k+1}
\end{align*}
with the last step following by induction.    Thus we have that
\[
\frac{d F_{m,k}(z,u)}{du}|_{u=1} = m T(z)^{m-1} (T(z) + T(z)^2 + \ldots + T(z)^{k+1}) \frac{z}{(1-T(z)) T(z)}.
\]
The leading term in the asymptotic expansion at $z= 1/e$ is $(k+1)m (e \sqrt{2} \sqrt{1 - ez})^{-1}$ again using \cite[Proposition 1]{FO90}, so applying \cite[Theorem 1]{FO90} and comparing with \eqref{eq:asymptoticfmn} gives the result.  \qed
\end{proof}

\begin{remark}
These results require that $k$ and $m$ be fixed as $n$ grows.  When analyzing iterating to supersingular $j$-invariants, this assumption does not hold: if there are $n$ $j$-invariants then on the order of $\sqrt{n}$ of them are supersingular.  
A more careful analysis should give that the asymptotic of Proposition~\ref{prop:expectedfixedpoints} continues to hold as long as $k$ and $m$ grow ``slowly'' compared to $n$.  Given the conclusion that iteration is not helpful for finding supersingular $j$-invariants, we do not pursue this.
\end{remark}

 In light of this analysis of functions with many fixed points, the iterative methods investigated in Section~\ref{ss:iterationrandom} behave exactly like random functions with the correct number of fixed points.

\section{Modular polynomials and curves isogenous to their conjugates}\label{sec:modular}

\subsection{Overview}

As described in the introduction and Section~\ref{sec:sec2}, Br\"oker's method is limited by
the degree of the Hilbert polynomials, upon which the runtime depends.
However, taking small-degree Hilbert polynomials leads to curves with small
endomorphisms (a vulnerability).
In this section,
we consider using polynomials 
whose roots correspond to curves with endomorphisms of exponentially large degree.  The hope is, at least, to demonstrate a hard curve.  The process we will describe does not, na\"ively, appear likely to generate curves in a uniformly random manner, although perhaps it can be adapated.

If $n$ is a positive integer coprime to $p$,
then the classical modular polynomial $\Phi_n(x,y) \in \ZZ[x,y]$
is defined as follows.
For any elliptic curve $E$,
let $S_{E,n} = \{ C \subseteq E[n] : C \text{ cyclic}, \#C = n \}$.
There are $\psi(n)$ elements of $S_{E,n}$, where $\psi$ is the Dedekind psi function
(recall that $\psi(\ell^k) = (\ell+1)\ell^{k-1}$ for $\ell$ prime, and
$\psi$ is multiplicative; in particular, \(\psi(n) > n\)).
Write $E/C$ for the codomain of a separable $n$-isogeny from $E$ with kernel $C$. Then
\[
   \Phi_n( j(E), y ) =  \prod_{C \in S_{E,n}} (y - j( E/C )).
\]
In other words, $\Phi_n(x,y) = 0$ if and only if $x$ and $y$ are
$j$-invariants related by a cyclic $n$-isogeny.
This remains the case over any field.  (See \cite[Chapter 5]{Lang} for background.)

Now, consider the roots in $\FF_{p^2}$ of the univariate polynomial
$\Phi_n(x,x^p)$.
These roots are the $j$-invariants of curves with
cyclic $n$-isogenies to their conjugates
(with root multiplicities equal to the number of distinct $n$-isogeny
kernels).
Let $j$ be such a root and let $E$ over $\FF_{p^2}$ be an elliptic curve with $j(E) = j$.
Denote by $E^{(p)}$ the Galois conjugate with respect to $\FF_{p^2}/\FF_{p}$.
One can compose the cyclic $n$ isogeny $\phi : E \to E^{(p)}$ with the inseparable Frobenius map $\pi : E^{(p)} \to E^{(p^2)} = E$.
When $\phi^{(p)} = \pm \hat{\phi}$, which is the general case, then~\cite[Proposition~2]{chenu2021Higher} shows that (possibly by taking a quadratic twist) we obtain an isogeny $\mu = \pi \circ \phi$ that satisfies $\mu^2 = -np$.
An alternative explanation for this is given in the proof of \cite[Lemma 6]{CGL}, where it is shown that if $n$ is small compared to $p$ and $E$ is supersingular then $\mu^2 = -np$.
Either way, it follows that the class group is of $\QQ(\sqrt{-np})$ acts on a large subset of the supersingular roots of $\Phi_n(x,x^p)$.


There is no particular reason why these curves should also have 
small-degree non-integer endomorphisms.

The collection of \emph{supersingular} curves with an $n$-isogeny to the
conjugate has been studied
\cite{Arpin,Supersingularland,chenu2021Higher}, and plays a role in the
security of the path-finding problem \cite{eisentrager2020computing}.
Since the class group of $\QQ(\sqrt{-np})$ acts on (a large subset of) this set, these curves form CSIDH-like graphs which could be used for
cryptographic purposes~\cite{chenu2021Higher}.
Thus,
a construction for random supersingular curves involving \(\Phi_n(x,x^p)\)
may lead to a means of sampling from these CSIDH-like graphs.
As in the CSIDH setting,
there are subexponential quantum algorithms to solve the \emph{vectorization} or \emph{class group action} problem (see \cite[Section 9.1]{ArpinEtAl21}, \cite{chenu2021Higher} and \cite{Wesolowski21Orientations}).  Thus, if there is a curve of known endomorphism ring in this set (see for example \cite{chenxue}),
one may be able to solve the fundamental isogeny problems (path-finding
and endomorphism ring computation) in quantum subexponential time.
This is still far from polynomial and may be considered secure for some
applications.

For \(p > n\), the polynomial \(\Phi_n(x,x^p)\) has degree \(\psi(n)p\),
which is  exponential with respect to \(\log p\).
While this polynomial is quite sparse, especially when \(p \gg n\),
we cannot compute its roots efficiently.
The idea is to reduce that degree, and make computations manageable,
by instead computing roots of the factor(s)
\[
    f_{n,m,p}(x) := \gcd( \Phi_n(x,x^p), \Phi_m(x,x^p) )
\]
for some auxiliary \(m\),
without explicitly computing 
\(\Phi_n(x,x^p)\) or \(\Phi_m(x,x^p)\).

The proposed approach for constructing supersingular curves is then:
\begin{enumerate}
    \item Choose $n$ and $m$.
    \item Compute one or more roots of  $f_{n,m,p}(x)$ in \(\FF_{p^2}\).
        There are \(O(nm)\) of these roots, and
        we can compute them in polynomial time with respect to \(n\), \(m\), and \(\log p\)
        (see~\S\ref{sec:fnmp-roots}).
    \item Test each root to see if it is a supersingular
        \(j\)-invariant,
        using e.g.~Sutherland's supersingularity
        test~\cite{2012/Sutherland};
        we give heuristics for this step in~\S\ref{sec:root-heuristics}.
\end{enumerate}

This method produces curves known to have endomorphisms of degree $nm, np$ and $mp$.
Since we wish to avoid endomorphisms of small degree,
the presence of the degree-\(mn\) endomorphism means that we should take at least one of $n$ and $m$ to be exponentially large.  
Nevertheless, it is plausible that the information about the endomorphism leaked from the process of construction is not enough to allow us to compute $\End(E)$ efficiently (i.e., in polynomial time).

\subsection{Computing roots of $f_{n,m,p}$}
\label{sec:fnmp-roots}

We want to compute 
roots 
of \(f_{n,m,p}(x) = \gcd(\Phi_n(x,x^p),\Phi_m(x,x^p))\)
in \(\FF_{p^2}\).
Note that simply computing
$\Phi_m(x,x^p)$ and $\Phi_n(x,x^p)$ in \(\FF_p[x]\), computing their \(\gcd\),
and finding its roots 
is exponential in \(\log p\),
because 
\(\deg\Phi_m(x,x^p) > mp\) and \(\Phi_n(x,x^p) > np\);
these polynomials are sparse for large \(p\),
but generic \(\gcd\) computations
(which are quasilinear in the maximum of the degrees of the inputs~\cite{Moller})
cannot take advantage of this.

Algorithm~\ref{alg:fnmp-roots} computes
all of the \(\FF_{p^2}\)-roots\footnote{%
    Algorithm~\ref{alg:fnmp-roots} ignores root multiplicities, but can be easily modified
    to take them into account if required.
}
of \(f_{n,m,p}(x)\)
in polynomial time with respect to \(m\), \(n\), and \(\log p\).
The key to its polynomial runtime in \(\log p\) 
is that the polynomials \(F_m\) and \(F_n\) constructed in
Lines~\ref{alg:fnmp-roots:F_m}
and~\ref{alg:fnmp-roots:F_n}
satisfy (by definition)
\[
    \Phi_m(j,j^p) = \Phi_n(j,j^p) = 0
    \iff
    F_m(j_0,j_1) = F_n(j_0,j_1) = 0
\]
for all \(j = j_0 + j_1\sqrt{\delta}\) in \(\FF_{p^2} =
\FF_p(\sqrt{\delta})\),
and it is much easier to solve 
the bivariate system \(F_m(X_0,X_1) = F_n(X_0,X_1) = 0\)
than it is to compute \(\gcd(\Phi_m(x,x^p),\Phi_n(x,x^p))\)
when \(p\) is large.

\begin{algorithm}
    \caption{%
        Compute the set of roots of
        \(f_{n,m,p}(x) = \gcd(\Phi_n(x,x^p),\Phi_m(x,x^p))\)
        in \(\FF_{p^2}\).
    }
    \label{alg:fnmp-roots}
    \KwIn{\(m\), \(n\), \(p\)}
    \KwOut{The set of roots of \(f_{n,m,p}(x)\) in \(\FF_{p^2}\)}
    Compute \(\Phi_m(X,Y)\) and \(\Phi_n(X,Y)\) in \(\FF_p[X,Y]\)
    \tcp*{Using e.g.~the algorithm of~\cite{BLS}}
    \label{alg:fnmp-roots:Phi}
    Compute a nonsquare \(\delta\) in \(\FF_p\),
    and its square root \(\sqrt{\delta}\) in \(\FF_{p^2}\)
    \;
    \(F_m \gets \Phi_m(X_0+\sqrt{\delta}X_1,X_0-\sqrt{\delta}X_1)\)
    in \(\FF_p[X_0,X_1]\)
    \tcp*{\(F_m\in\FF_p[X_0,X_1]\) because \(\Phi_m\) is symmetric}
    \label{alg:fnmp-roots:F_m}
    \(F_n \gets \Phi_n(X_0+\sqrt{\delta}X_1,X_0-\sqrt{\delta}X_1)\)
    in \(\FF_p[X_0,X_1]\)
    \tcp*{\(F_n \in \FF_p[X_0,X_1]\) because \(\Phi_n\) is symmetric}
    \label{alg:fnmp-roots:F_n}
    \(R \gets \) \Resultant{\(F_m\), \(F_n\); \(X_0\)}
    \tcp*{Bivariate resultant \(\mathrm{Res}_{X_0}(F_m,F_n)\) in \(\FF_p[X_1]\)}
    \(\mathcal{J}_1 \gets \) \Roots{\(R\), \(\FF_p\)}
    \;
    \(\mathcal{S} \gets \emptyset\)
    \;
    \For{\(j_1\in \mathcal{J}_1\)}{
        \(G \gets \) \GCD{\(F_m(X_0,j_1)\), \(F_n(X_0,j_1)\)}
        \;
        \(\mathcal{J}_0 \gets \) \Roots{\(G\), \(\FF_p\)}
        \;
        \(
            \mathcal{S}
            \gets
            \mathcal{S} \cup \{j_0 + j_1\sqrt{\delta} : j_0 \in \mathcal{J}_0\}
        \)
    }
    \Return{\(\mathcal{S}\)}
\end{algorithm}

As has already been noted, for security in applications, at least one of $n$ and $m$ must be exponentially large. 
But if \(n\) (or \(m\)) is super-polynomially large with respect to \(\log p\),
then Algorithm~\ref{alg:fnmp-roots} requires super-polynomial time and
space, since it must work explicitly with the polynomials \(\Phi_n\).
Hence a natural question is whether we can do better than Algorithm~\ref{alg:fnmp-roots} when one (or both) of $n$ and $m$ is large.
This is an open question.
If $m$ is small and $n$ is very large then a ``dream'' approach would be to compute $F_m$ using the classical algorithm and then somehow compute $\text{Resultant}( F_m , F_n ;  X_0 )$ directly by some form of ``square-and-multiply'' approach without explicitly computing $F_n$.

\subsection{Supersingular roots of $f_{n,m,p}$}
\label{sec:root-heuristics}

Now we consider the question of how many of the roots of $f_{n,m,p}(x)$
might be supersingular \(j\)-invariants.
The individual polynomials $\Phi_n(x,x^p)$ should
be expected to have overwhelmingly ordinary roots,
but there are some heuristic reasons to expect \(f_{n,m,p}(x)\)
to have a higher proportion of supersingular roots and we give some evidence for this in Section~\ref{sec:experiments}.

A heuristic lower bound on the number of supersingular roots of 
$f_{n,m,p}(x)$
can be obtained as follows.
There are \(\approx p/12\) supersingular curves over \(\FF_{p^2}\),
and \(\approx \sqrt{np}\) of them have an \(n\)-isogeny to their conjugate
(combining \cite[Theorem 2]{chenu2021Higher} and heuristic average class group estimates).
Hence, we can postulate that a ``random'' supersingular curve has probability $\sqrt{12n/p}$ of having an \(n\)-isogeny to its conjugate.
Applying this to the \(\approx \sqrt{mp}\) supersingular curves that are roots of \(\Phi_m(x,x^p)\), we conclude that there should be \( \approx \sqrt{mn} \) supersingular roots of $f_{n,m,p}(x)$.

To summarize, if we expect that the properties of having an $n$-isogeny and having an $m$-isogeny to the conjugate are in an appropriate sense ``independent,'' 
then one might expect the supersingular portion of $f_{n,m,p}(x)$ to have degree $\approx \sqrt{nm}$. 
 Note that the resultant in line 5 of Algorithm~\ref{alg:fnmp-roots} has degree $O(mn)$, so
there are \(O(mn)\) roots of \(f_{n,m,p}(x)\) in \(\FF_{p^2}\)
(to see this, apply Bézout's theorem to the polynomials \(F_m\) and
\(F_n\) in Algorithm~\ref{alg:fnmp-roots}). 

Given the degree estimate just described, one might consider taking the gcd of three different modular polynomials.  This will almost certainly have a smaller degree:  continuing the heuristic argument above would lead to $O(\sqrt{nmr/p})$ supersingular roots for the gcd of $\Phi_n(x,x^p)$, $\Phi_m(x,x^p)$ and $\Phi_r(x,x^p)$.   With such an estimate, one might consider taking $n \sim m \sim \sqrt{p}$ and $r$ polynomial in $\log p$.   One might expect the 3-way gcd to have supersingular roots, provided it is not $1$, by the same heuristics as above.

\begin{remark}
    If $E$ has an $n$-isogeny and an $m$-isogeny to its
    conjugate $E^{(p)}$, then it also has an $nm$-endomorphism to itself.
    When $p$ is inert in $\QQ(\sqrt{-nm})$,
    some such $E$ will be reductions modulo $p$
    of curves over $\overline{\mathbb{Q}}$ with CM by $\ZZ[\sqrt{-nm}]$,
    specifically those where the reduction of the $nm$-endomorphism factors
    through the conjugate.
    In this case, we can expect a nontrivial gcd between $f_{n,m,p}(x)$
    and the Hilbert class polynomial for $\mathbb{Q}(\sqrt{-nm})$.
\end{remark}

\begin{remark}
If one desired uniformly randomly generated supersingular curves, one might consider using randomly generated $n$ and $m$ of a certain size.  It is unclear what distribution of $n$ and $m$ would lead to a uniformly random distribution of supersingular curves, if any.
\end{remark}

\subsection{Experimental evidence}\label{sec:experiments}

\begin{figure}[h!]
\begin{center}
    \includegraphics[width=0.45\textwidth]{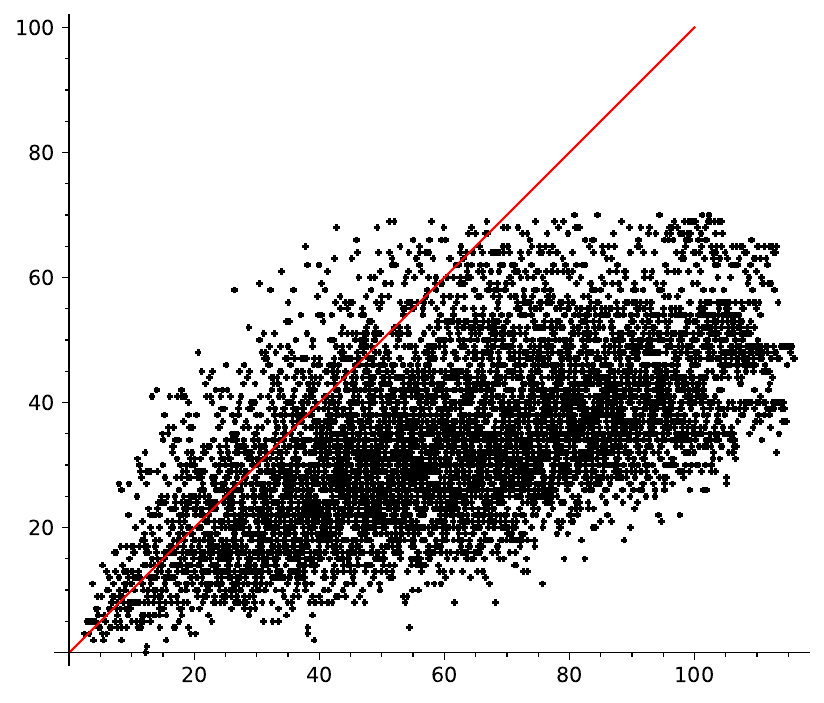} 
\end{center}
\caption{Scatterplot of $f_{n,m,983}$:  $x$-axis is $\sqrt{nm}$, $y$-axis is the number of $\FF_{p^2}$ roots.  The line $y=x$ is shown for reference.  There are a total of 7046 data points.}
\label{fig:degratio}
\end{figure}

\begin{figure}[h!]
\begin{center}
    \includegraphics[width=0.45\textwidth]{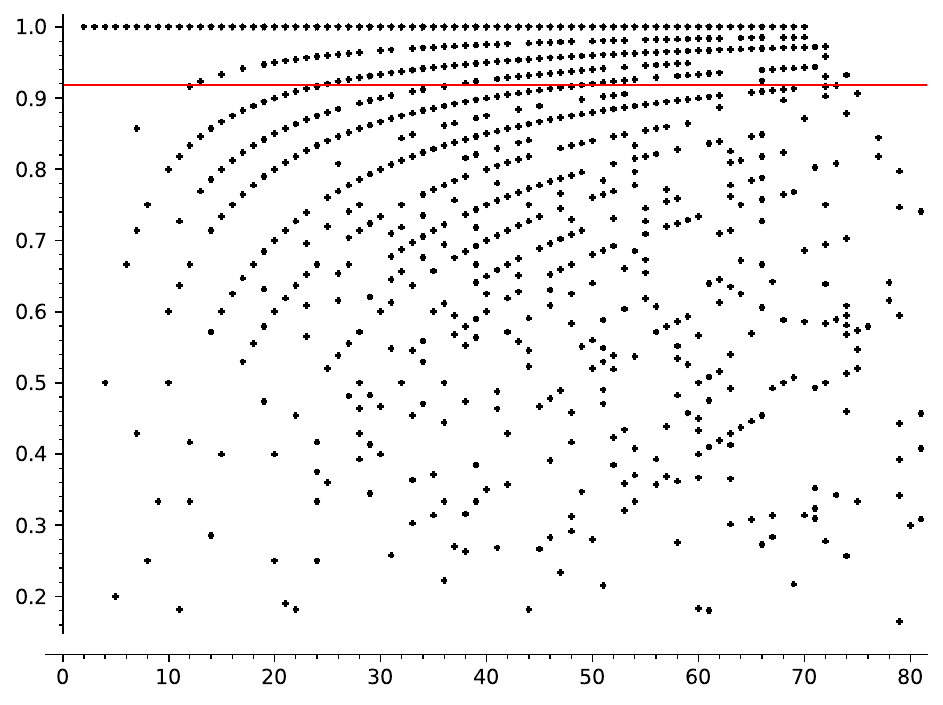}
    \includegraphics[width=0.45\textwidth]{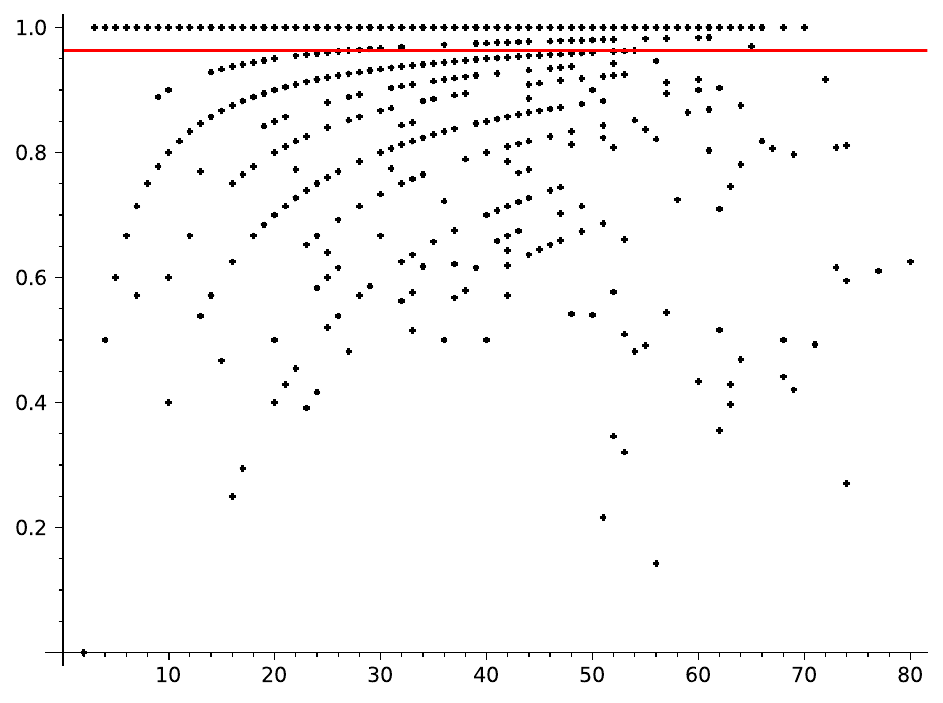}
\end{center}
\caption{Scatterplot of $f_{n,m,983}$:  $x$-axis is degree, $y$-axis is the ratio of supersingular roots to all $\FF_{p^2}$ roots.  The visible hyperbolas correspond to the existence of $0$, $1$, $2$ etc. non-supersingular roots.  At left, 4286 pairs $(n,m)$ which are coprime; at right, 2760 pairs which are not coprime.  The red line indicates the average ratio across all pairs.}
\label{fig:ssratio}
\end{figure}

\begin{figure}[h!]
\begin{center}
    \includegraphics[width=0.45\textwidth]{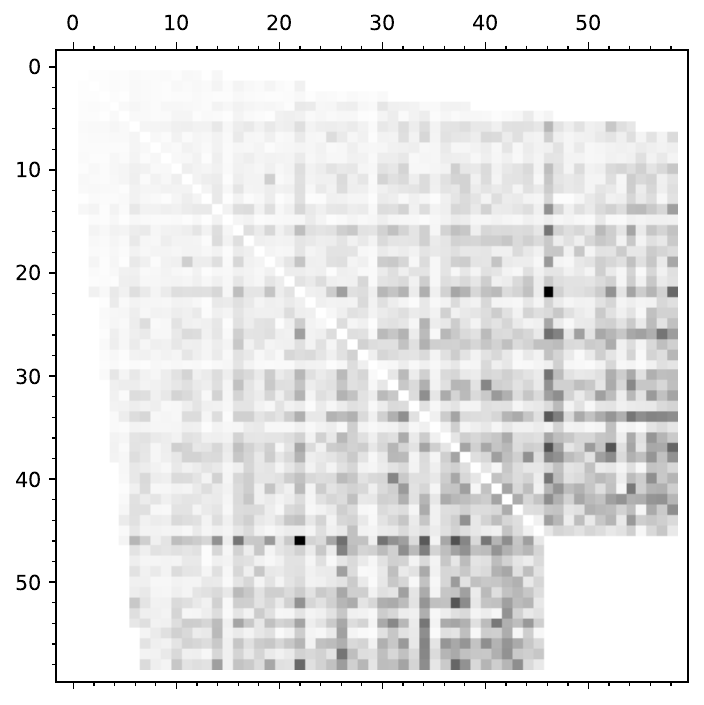}
    \includegraphics[width=0.45\textwidth]{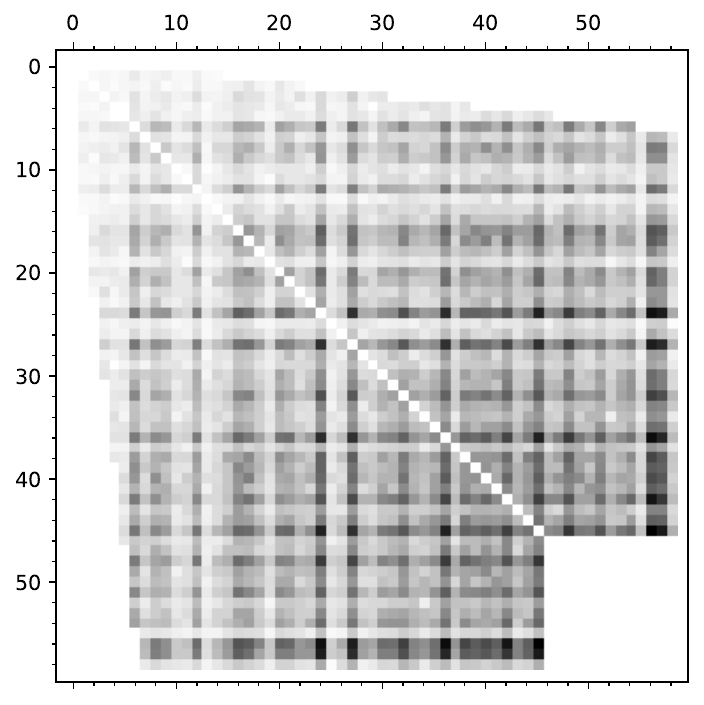}
\end{center}
\caption{Plot of the number of $\FF_{p^2}$ roots (left) and supersingular $\FF_{p^2}$ roots (right) of $f_{n,m,983}$ as a function of $n$ and $m$ ($x$ and $y$ axes).  Dark = more (maximum = 238 $\FF_{p^2}$ roots at left; 70 supersingular roots at right); light = fewer (minimum = 0); white = uncomputed.}
\label{fig:matplot}
\end{figure}

\begin{figure}[h!]
\begin{center}
    \includegraphics[width=0.45\textwidth]{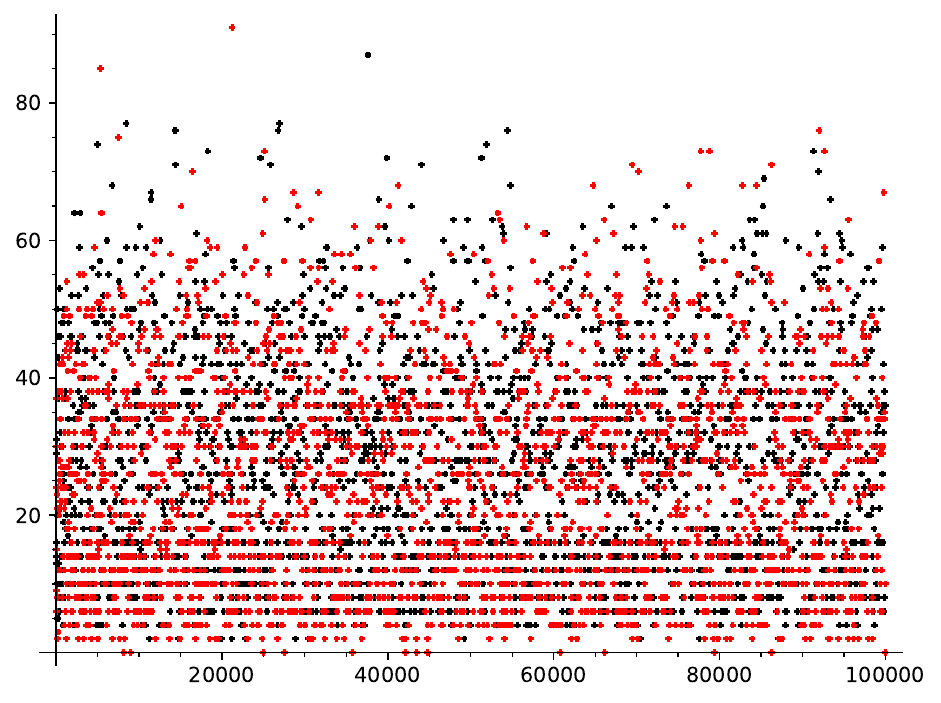}
    \includegraphics[width=0.45\textwidth]{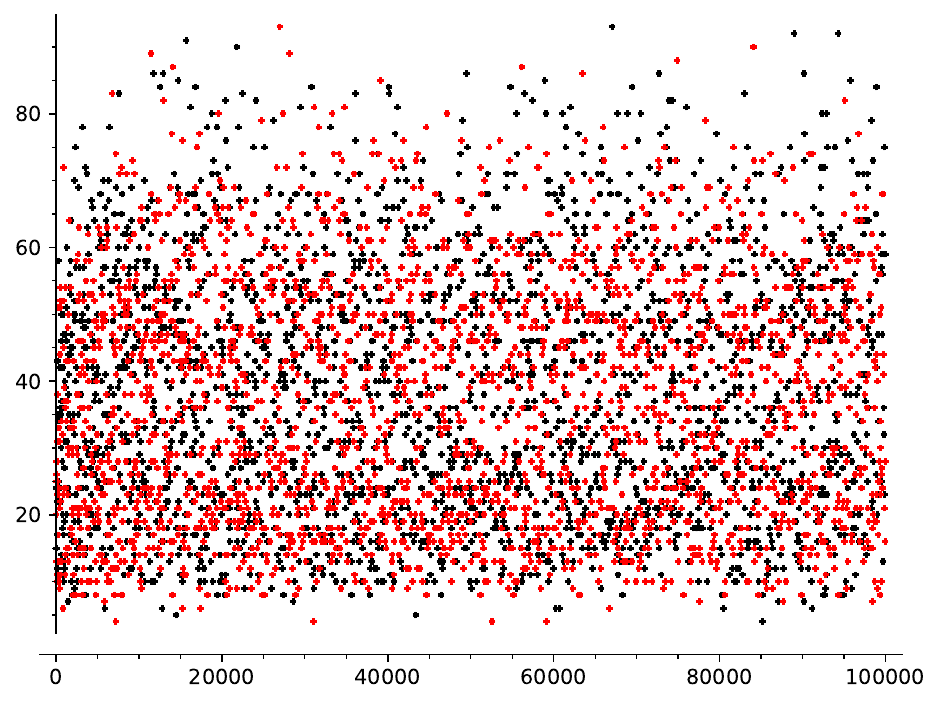}
\end{center}
\caption{Plots of $\deg f_{n,m,p}$ as a function of $p$.  At left, $n=8$, $m=12$.  At right, $n=8$, $m=13$.  When $p$ is inert in $\QQ(\sqrt{-nm})$, the plotted points are black.  When $p$ is split, the plotted points are red.  In both plots, there are 4808 points in total, representing a random selection of primes which are $3 \pmod 4$ in the given range.
}
\label{fig:pchanging}
\end{figure}

\begin{table}[]
    \centering
    \setlength{\tabcolsep}{0.5em}
\renewcommand{\arraystretch}{1.8}
    \begin{tabular}{c|c|c|c}
         data set & points & avg. $\frac{\text{supersingular roots}}{\text{$\FF_{p^2}$ roots}}$ & avg. $\frac{\text{$\FF_{p^2}$ roots}}{\sqrt{nm}}$  \\
         \hline
         all & 7046 & 0.9362 & 0.6251 \\
         coprime & 4286 & 0.9189 & 0.6555  \\
         not coprime & 2760 & 0.9632 &  0.5779
         \\
         $\left( \frac{-nm}{p} \right) = -1$ & 3572 & 0.9424 & 0.6373 \\
         $\left( \frac{-nm}{p} \right) = 1$ & 3474 & 0.9300 & 0.6126\\
    \end{tabular}
    \caption{Statistics for various subsets of the data set for $p=983$.  The first row (`all') contains all the polynomials $f_{n,m,983}$ collected as described in the beginning of the section.  The other rows give statistics for subsets of the data where $n$ and $m$ satisfy some criterion:  the row `coprime' (respectively `not coprime') refers to those data points where $\gcd(n,m)=1$ (respectively, $\gcd(n,m)\neq 1$), and the final two rows include points where $n$ and $m$ satisfy the indicated equality.}
    \label{tab:ss-data}
\end{table}


To test the heuristics of the previous section, the polynomials $f_{n,m,p}(x)$ were computed for a fixed prime $p = 983$ with pairs $(n,m)$ ranging over $2 \le n \le 45$, $n + 1 \le m < 8n$. Figure~\ref{fig:degratio} shows the degree of the $\FF_{p^2}$ part of the polynomial as compared with $\sqrt{nm}$.  Figure~\ref{fig:ssratio} shows the proportion of supersingular roots.  Table~\ref{tab:ss-data} gives the average values of these quantities for various subsets of the data, including where $(n,m)$ is coprime or not, and where $p$ is inert or split in $\QQ(\sqrt{-nm})$. Figure~\ref{fig:matplot} gives a sense of how the number of $\FF_{p^2}$ roots of $f_{n,m,p}$ varies in an intricate manner as a function of $n$ and $m$ for fixed $p$.  

In addition, the polynomials $f_{n,m,p}(x)$ were computed for various fixed pairs $(n,m)$ with $p$ ranging over all primes less than $10^6$ with $p \equiv 3 \pmod{4}$. Figure~\ref{fig:pchanging} shows the degrees of $f_{8,12,p}$ and $f_{8,13,p}$ with respect to $p$.

In general, the data seems to support the following patterns: \begin{enumerate*}[label=(\roman*)] \item the degrees of the $\FF_{p^2}$ parts of $f_{n,m,p}$ may be similar to or slightly less than $\sqrt{mn}$, \item the proportion of supersingular roots among $\FF_{p^2}$ roots is often high, \item there is variation in the ratio of supersingular roots with $n$ and $m$, with slightly higher proportions found amongst $n$ and $m$ not coprime, and \item as $p$ varies, the degree of $f_{n,m,p}$ is relatively constant, but is dependent upon the coprimality, not just size, of $n$ and $m$.
\end{enumerate*}

To conclude, supersingular $j$-invariants seem to be a large proportion of the $\FF_{p^2}$ roots of $f_{n,m,p}$. Hence, if there were an efficient way to compute these roots for a wide range of values $(m,n)$ then this might give a solution to the problem of hashing to a supersingular curve.

\section{Constructing supersingular curves using constraints on their torsion
}\label{sec:reverse}

A supersingular curve is characterized by the number of points over any extension. Provided a curve, Schoof's algorithm~\cite{schoof1985elliptic} computes the trace. When hashing into supersingular graphs, we know the trace and we want to find a curve. Thus, one may try to use Schoof's algorithm ``backwards,'' by setting up a system of equations restricting the trace (or, more directly, the field of definition of torsion points), and looking for solutions.  This method may lead to a way to generate supersingular curves uniformly randomly, since some such systems have all supersingular curves as roots. 

\subsection{A system of equations}
\label{sec:rev-system}

To introduce the approach, let us first discuss the case when $p$ is a prime of the form $p+1=\prod_i\ell_i$, where $\ell_i$ are small distinct odd primes. For such $p$, the approach could proceed as follows.
Let $a$ be some parameter for the curve, like the $j$-invariant or the Montgomery coefficient. For every~$i$, write $\Psi_{\ell_i}(x_{\ell_i},a)$ for the division polynomial of order $\ell_i$ of the curve parameterized by $a$. These polynomials can be efficiently computed. 
Consider the system 
\begin{equation}
    \label{eqn:rev-sch-1}
\left\{\begin{array}{lcl}
\Psi_{\ell_i}(x_{\ell_i},a)&=&0 \qquad \forall\ \ell_i|p+1 \\
x_{\ell_i}^{p^2}-x_{\ell_i}&=&0 \qquad \forall\ \ell_i|p+1,
\end{array}\right.
\end{equation}
with variables $x_{\ell_i}$ and $a$.
The equations of this system force the $\ell_i$-torsion points 
of the curve with parameter $a$ to be defined over $\mathbb{F}_{p^2}$ for all~$i$. %
Therefore the $p+1$ torsion is also defined over $\mathbb{F}_{p^2}$, which implies that any curve with parameter $a$ being a solution of \eqref{eqn:rev-sch-1} is supersingular.
Taking the resultant of all polynomials in the system with respect to all variables but~$a$ gives a polynomial whose roots are all parameters~$a$ that correspond to supersingular curves. 

More generally when $p+1$ is not smooth, one can fix a set of small primes or prime powers $\ell_i$ such that their product is above the Hasse bound, and replace the equations $x_{\ell_i}^{p^2}-x_{\ell_i}=0$ in \eqref{eqn:rev-sch-1} by alternative equations forcing the endomorphisms $$\pi^2+[p-1]\pi+p^2,$$
where $\pi$ denotes the Frobenius endomorphism, on the curve with parameter $a$ to act trivially on the $\ell_i$ torsion.

For primes of the form $p+1=f\ell_2^{e_2}\ell_3^{e_3}$ where $f,\ell_2,\ell_3$ are small integers (as used in the SIDH key exchange ~\cite{jao2011towards}) one can replace a single equation
$\Psi_{\ell_i}(x_{\ell_i},a) = 0$ in \eqref{eqn:rev-sch-1} by a polynomial system in the variables $x_{ij}$ and $a$
\begin{equation} \label{eq: polys Fp2}
\left\{\begin{array}{lcl}
\Psi_{\ell_i}(x_{i1},a) & = & 0 \\
\mbox{}[\ell_i]_a(x_{i(j+1)},-) & = & (x_{ij},-) \quad \mbox{ for } 1\leq j \leq e_i-1,\\
\end{array}\right.
\end{equation}
where $[\ell_i]_a$ are ``$x$-only'' multiplication-by-$\ell$ polynomials on the curve of parameter $a$. For any solution to this system,  $x_{ij}$ is the $x$-coordinate of a point $(x_{ij},-)$ of order $\ell_i^j$ on the curve with parameter $a$. Note that the equations $[\ell_i]_a(x_{i(j+1)},-)  =  (x_{ij},-)$ are of degree roughly $\ell_i^2$ and $\Psi_{\ell_i}(x_{i1},a)$ is of degree $(\ell_i^2-1)/2$.

\medskip
As with other approaches involving large polynomial systems or large degree equations, the cost and optimal strategy to solve these systems are not obvious. We observe that the polynomial system~(\ref{eq: polys Fp2}) 
contains equations in $e_1+e_2+1$ variables of degree roughly $\ell_i^2$ and $\ell_i$ together with the equation translating the fact that the torsion points lie in $\mathbb{F}_{p^2}$
of degree $p^2$.
Yet, compared to generic polynomial systems of the same degree and with an equal number of variables, the given polynomial systems have only a few mixed monomial terms. Further, they exhibit a certain block structure. Instead of using generic algorithms such as Gröbner basis computations, taking the full monomial structure into account might help to solve the polynomial systems faster. This might be feasible using algorithms such as Rojas' algorithm for sparse polynomial systems~\cite{rojas1999solving}. However, further research is needed to draw conclusions about the concrete speedup that can be achieved using this additional structure and to assess the cost of solving the polynomial systems given in this section.

Unlike the Hasse polynomial of Section~\ref{sec:Iterating} and arguably the function $f_{n,m,p}$ of Section~\ref{sec:modular}, the polynomial system to be solved in this section can certainly be stored in polynomial space.

\subsection{Variants}

\subsubsection{Reducing the number of solutions:} Instead of computing a random solution to the polynomial systems described in the previous section and thus a random curve with the correct number of points, some applications require computing
only one curve with unknown endomorphism ring. To achieve this, one could add additional equations to the systems~(\ref{eq: polys Fp2})
to reduce the number of expected solutions -- potentially all the way to 1, when solving the system would mean to select a single curve. 

\medskip
One approach could be to restrict the $x$-coordinate of torsion points to random cosets of multiplicative subgroups, namely replacing
$x_{\ell_i}^{p^2}-x_{\ell_i}=0$ for some $i$ by
$$(\mu_ix_{\ell_i})^{r_i}-1=0$$
for suitable $r_i$ dividing $p^2-1$, and random $\mu_i$ in $\mathbb{F}_{p^2}$. This will decrease the degrees of equations in the system, as well as the number of solutions. If one does not restrict the field equations for all $i$, one may want to choose some $i$ uniformly at random.

\medskip
Assuming that the solutions to the system~(\ref{eq: polys Fp2}) are ``randomly'' distributed among all cosets of the multiplicative subgroup, the expected number of solutions to the system is reduced by the number of such cosets. If one of the remaining solutions is chosen uniformly at random and if the cosets for different $i$ were chosen uniformly at random, then the supersingular elliptic curve corresponding to the final solution is a random supersingular elliptic curve.
 One could consider various versions of this, leaving more or fewer solutions.

 \subsubsection{Hybrid version:} Another variant is to drop some equations in the polynomial system~(\ref{eq: polys Fp2}). The resulting system has then more solutions. Each solution to the resulting system leads to a curve with a number of points $N$ with trace not fixed modulo the Hasse bound. That is, the curve generated might be of order $N$ different from 
 the order $(p + 1)^2$ we would like to find. Hereby, the number of equations dropped from the system~(\ref{eq: polys Fp2}) controls the size of $\gcd(N,(p+1)^2)$.
Thus, to compute a supersingular elliptic curve one may want to proceed as follows. One generates a system with fewer equations and keeps computing random solutions until the resulting curve has the correct order. We leave it for future research to examine how much easier it is to solve the resulting systems compared to~(\ref{eq: polys Fp2}).

\section{Genus 2 Walks}\label{sec:genus_2}

In this section, we explore several approaches to sample a uniformly random supersingular elliptic curve based on the following general idea: start with a known supersingular elliptic curve $E_0/\FF_q$, glue it to itself to construct a genus-2 Jacobian $A \cong \Jac(C)$ explicitly isogenous to $E_0^2$, and then connect $A$ with a new random-looking elliptic product using Richelot isogenies, or through geometric inspection of the Jacobian (via its Kummer surface).
One might hope that these genus-2 operations will ``hide'' obvious isogenies between the elliptic curves involved, but we will explain a number of issues with this approach at the end of Section~\ref{Sec:RandomWalks}.

Let $A$ be a principally polarised abelian surface (PPAS)
over a finite field \(\FF_q\) of characteristic \(p > 2\).
The correct generalisation of the notion of supersingularity to genus 2
is to say that $A$ is supersingular if and only if the Newton polygon of its
Weil polynomial has all its slopes equal to $1/2$;
this is the case if and only if the $p$-torsion $A[p]$ is isomorphic 
(as a $\mathrm{BT}_1$ group scheme)
to either $I_{2,1}$ or $I_{1,1} \oplus I_{1,1}$, where $I_{1,1}$ is the
$p$-torsion group scheme of a supersingular elliptic curve
(see Pries~\cite{pries2008} for further detail).
In the latter case, we say $A$ is principally polarized \emph{superspecial} abelian surface (PPSSAS).

Every PPAS $A$ is isomorphic (as a principally polarized abelian variety)
to either the Jacobian $\Jac(C)$ of some genus-2 curve $C$, or the product $E_1 \times E_2$ of
two elliptic curves (which are both supersingular if $A$ is superspecial). 
Oort~\cite{oort1975} has shown that every superspecial abelian surface
is isomorphic \emph{as an unpolarized abelian variety} to a product of
supersingular elliptic curves, and that every supersingular abelian
surface is at least isogenous to a product of supersingular elliptic
curves (if the abelian surface is supersingular but not superspecial,
then the isogeny is inseparable).

We can construct a superspecial Jacobian $A$ isogenous to
a product of supersingular elliptic curves $E_1$ and $E_2$ by
\emph{gluing} them along their $2$-torsion, say. This corresponds to a
Richelot isogeny~\cite{richelot1837} $E_1\times E_1 \to (E_1 \times E_2)/G \cong A$,
where $G \leq (E_1 \times E_2)[2]$ is the graph of an isomorphism of
group schemes $\psi: E_1[2] \to E_2[2]$ that is an anti-isometry with
respect to the $2$-Weil pairing (see \cite{Kani1997});
the resulting $A$ is always a Jacobian.
(We can also glue along the $\ell$-torsion for $\ell > 2$,
and there is an analogous inseparable construction in~\cite{oort1975}
for gluing along the $p$-divisible group schemes $E_i[\operatorname{Fr}_p]$,
but the case $\ell = 2$ is sufficient to illustrate our ideas. There is
no reason to suspect that $\ell > 2$ or $\ell = p$ will give better
results. The case $\ell=2$ also has the advantage of being completely explicit.)

\subsection{Random Walks}\label{Sec:RandomWalks}
Our first idea is simple: We begin with a supersingular elliptic curve and glue it to itself which induces an isogeny to an abelian surface. 
We then take a random walk on the isogeny graph of abelian surfaces. 
Finally, we find the closest reducible surface and return one of its supersingular elliptic factors.
The idea can be summarised in the following diagram:
\[E\times E\xrightarrow{\text{ glue }} A \xrightarrow{\text{ rand. walk }} A'\xrightarrow{\text{ unglue }} E'\times E''\]
The initial \(A\) is superspecial, and so superspeciality is preserved so long as the isogenies in the random walk are of degree prime to the characteristic. 
This means that we are walking in the superspecial graph.

A similar situation occurs in \cite{CostelloSmith2020}, where the
authors consider the supersingular isogeny problem in genus 2 and
higher. We will only sketch the outline of their arguments and will
refer interested readers to find details in their paper: In genus 2,
given two superspecial abelian surfaces $A$ and $A'$, the idea is to
reduce the problem of finding an isogeny $\phi: A \to A'$ to the problem
of finding a factored isogeny $\psi: E_1 \times E_2 \to E_1' \times
E_2'$ and (un)gluings $\pi:A \to E_1 \times E_2$ and $\pi':E_1' \times
E_2' \to A'$. Finding the isogenies $\pi$ and $\widehat{\pi'}$ is
essentially done by taking random walks of length $O(\log(p))$. Such a
walk encounters a product of elliptic curves with probability $O(1/p)$,
so after $O(p)$ many random walks we should have found the
required $\pi$ and $\widehat{\pi'}$.
(The heuristics of~\cite{CostelloSmith2020} are made more rigorous
in~\cite{FloritSmith2022}.)

Translating this to our setting, we see that random walks away from a
fixed superspecial abelian surface have no better expected runtime at
encountering a supersingular elliptic curve than simply searching for
one directly by randomly sampling $j$-invariants and testing if they correspond to supersingular elliptic curves.

Ultimately, for this approach to give any advantage over simply taking a random walk
in the elliptic supersingular graph,
we need the genus-2 walk to ``hide'' information about the relative
endomorphism rings of the starting and ending elliptic factors.
But as noted in \cite[Section~2]{ibukiyamakatsuraoort1986}, by fixing a supersingular elliptic curve over a finite field it is possible to parametrise the space of PPSSASs by positive-definite hermitian matrices which are elements of the matrix algebra $M_2(B_{p,\infty})$, where $B_{p,\infty}$ is the definite quaternion algebra that is ramified exactly at $p$ and $\infty$.
Furthermore, isogenies between PPSSASs can be represented by conjugation by matrices in the same matrix algebra.
Thus, knowledge of the random walk in the genus-2 graph may allow the construction of a matrix in  $M_2(B_{p,\infty})$ that can be used to construct a path between our base and final supersingular elliptic curves. 

Lastly, knowledge of the genus-2 walk may allow for the adversary to compute the endomorphism ring of the target surface, by computing the matrix that corresponds to the isogeny walk. 
The endomorphism ring of an elliptic product contains the endomorphism
ring of each factor as a direct summand, so
this information should allow an adversary to compute the endomorphism ring of the resulting (supersingular) elliptic curve.

\subsection{Constructing curves on the Kummer surface}\label{sec:ungluing_P3}

We saw above that random walks in the superspecial genus-2 graph give no real
advantage over random walks in the elliptic supersingular graph when
constructing new supersingular elliptic curves---and in addition, they
may reveal information about the endomorphism ring.
But we know that every superspecial abelian surface $A$ is isomorphic to
an elliptic product \emph{as an unpolarised abelian variety},
so why not go looking for a new supersingular elliptic curve directly in $A$?

From a computational point of view,
it is easier to work with curves on the \emph{Kummer surface},
which is the quotient of $A$ by the action of the involution $[-1]$.
The projective embeddings of the Kummer surface $A/\langle{\pm1}\rangle$ 
are easier to manage than those of the abelian surface $A$,
since they involve fewer equations and lower-dimensional ambient spaces;
but they also retain much of the information of $A$.

In this part, we consider the singular model $\Ksing$ in $\PP^3$ 
of the Kummer surface of an abelian surface $A$.
The model $\Ksing$ is defined by a single quartic equation
(see e.g.~\cite[Eq. 3.1.8]{CasselsFlynn}).
We write $\pi: A \to \Ksing = A/\langle\pm1\rangle$ for the degree-two quotient map; this map
is ramified precisely at the sixteen $2$-torsion points of $A$, and the
images of these points under $\pi$ are the singular points of $\Ksing$,
known as \emph{nodes}. We denote the set of nodes by $S \subset \Ksing$.

If $E \subset A$ is an elliptic curve,
then the restriction of $\pi$ to $E$ defines a double cover of curves
$\pi: E \to E' := \pi(E) \subset \Ksing$.
It follows from the Riemann--Hurwitz formula that 
$E'$ is either an elliptic curve or a genus-0 curve;
$E'$ is an elliptic curve if and only if $\pi$ is unramified along $E$; 
and $E'$ is a genus-$0$ curve if and only if $\pi$ is ramified at precisely $4$ points.

This observation provides two ideas for constructing a new supersingular
elliptic curve from a superspecial abelian surface $A$:
\begin{enumerate}
    \item Find an elliptic curve on $\Ksing$ that does not go through
        any of the nodes of $\Ksing$. 
    \item Find a genus-$0$ curve on $\Ksing$ that goes through precisely
        $4$ of the nodes of $\Ksing$.
\end{enumerate}
For both approaches, we consider the intersection of $\Ksing$ with a hyperplane $H$.

\textbf{Approach 1:}  
For any hyperplane $H \subset \PP^3$, 
the intersection $\Ksing \cap H$ is a plane quartic curve $C$.
If $C$ is non-singular then it is a genus-$3$ curve. If on the other
hand $C$ is singular and has precisely two nodes then its (geometric) genus is $3-2 = 1$.
Hence, it is possible to obtain such genus-$1$ curves by constructing
hyperplanes that contain precisely two of the nodes of $\Ksing$.
Each pair of nodes determines a one-parameter family of
hyperplanes passing through them, and imposing singularity of the
intersection $C$ at the nodes gives simple algebraic conditions on the
parameter that let us choose ``good'' hyperplanes.
(If required, one may define a birational map from $C$ to an elliptic curve in Weierstrass form.)
There is an important caveat here: even if $C$ has genus 1, it may not
be the image of an elliptic curve in $A$.

In our experiments, we took $\Ksing$ to be the Kummer surface
of the Jacobian of the superspecial curve $y^2 = x^6 - x$ over~$\FF_{p^2}$
with $p \equiv 4 \pmod{5}$.
We note that this Jacobian is \emph{not} Richelot-isogenous to any
elliptic product (see e.g.~\cite[\S4.15]{FloritSmith2022Atlas}),
so we can be confident that any elliptic curves we find are not
connected with some gluing along \(2\)-torsion.
Unfortunately, none of the elliptic curves we found using this approach were supersingular.
We discuss reasons for this in~\S\ref{sec:why-ordinary} below.

\textbf{Approach 2:}
This approach is doomed to fail:
it is impossible to construct a hyperplane $H$ passing through
precisely $4$ of the nodes of $\Ksing$.
Any three of the singular points in $\Ksing$ already define a
hyperplane $H$, and it turns out that this hyperplane must pass through
exactly $6$ of the nodes.
These hyperplanes, known as the \emph{tropes} of the Kummer,
are classical objects of study; there are sixteen of them,
and the incidence structure formed by the intersections of tropes and
nodes is a $(16,6)$-configuration~\cite[\S26]{Hudson}.

If $H$ is a trope, then it is tangent to $\Ksing$.
The intersection is a smooth conic, taken twice,
and the preimage of this conic in $A$ is isomorphic to the genus-2 curve
generating $A$ as a Jacobian; its Weierstrass points are the ramified
points above the six nodes (see~\cite[\S3.7]{CasselsFlynn} for further
details, including the explicit recovery of the genus-2 curve).
This curve may degenerate to a union of two elliptic curves joined at one point,
but then $A$ is an elliptic product itself,
and these two elliptic curves are isomorphic to the factors---so
we cannot obtain any new supersingular elliptic curves in this way.

\subsection{Genus-$5$ curves on the desingularised Kummer}\label{sec:ungluing_P5}

We can find more elliptic curves by computing the 
desingularization $\phi: \Ksmooth \to \Ksing$ of the
Kummer surface, which yields a smooth model $\Ksmooth$ in $\PP^5$
(see~\cite[Chapter~16]{CasselsFlynn} for more details).
Concretely, 
let $Y:y^2=\prod_{i=0}^{5} (x-a_i)$ be a hyperelliptic curve.
Then $\Ksmooth = V(\Omega_0,\Omega_1,\Omega_2) \subset \PP^5$, where 
\[
\Omega_0: \sum_{i=0}^{5} X_i^2=0, \quad \Omega_1: \sum_{i=0}^{5} a_i X_i^2=0,\quad \Omega_2: \sum_{i=0}^{5} a_i^2 X_i^2=0
    \,,
\]
is a smooth model of the Kummer surface of the Jacobian variety of $Y$ 
(see Klein \cite{klein1870theorie}, and the survey articles by Dolgachev \cite{dolgachev2019kummer} and Edge \cite{edge1967new}). As an intersection of three quadrics in $\PP^4$, the intersection of $\Ksmooth$ with a hyperplane is a non-hyperelliptic genus-$5$ curve $C$.
We first explain how to construct different elliptic curves that arise as quotients of the curve $C$, and later explore an alternative path where we choose hyperplanes in such a way that the curve $C$ is singular and its irreducible components are elliptic curves. 

\subsubsection{Elliptic curves as quotients}
\begin{figure}
    \centering
    \begin{forest}
    for tree={inner sep=2pt,l=10pt,l sep=10pt}
    [$\Ksmooth$
        [$C_0$
            [$E_{0,1}$][$E_{0,2}$][$E_{0,3}$][$E_{0,4}$][$E_{0,5}$]
        ]
        [$C_1$
            [.][.][.][.][.]
        ]
        [$C_2$
            [.][.][.][.][.]
        ]
        [$C_3$
            [.][.][.][.][.]
        ]
        [$C_4$
            [.][.][.][.][.]
        ]
        [$C_5$
            [.][.][.][.][.]
        ]
    ]
    \end{forest}
    \caption{Elliptic curves $E_{i,j}$ contained in the Kummer surface $\Ksmooth$}
    \label{fig:construction}
\end{figure}
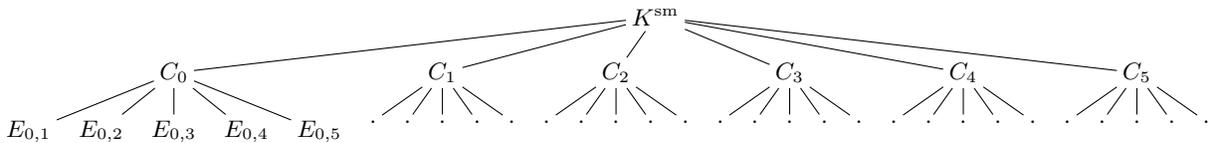

The intersection of the variety $\Ksmooth$ with a hyperplane defined by $X_i =0$ for some $i \in \{0,\dots,5\}$ yields a non-hyperelliptic genus-$5$ curve $C_i$. We are interested in certain elliptic curves $E_{i,j}$ with $j \in \{0,1,2,3,4,5\} \setminus \{i\}$ that arise as quotients of the curve $C_i$. This situation is  also studied by Stoll in \cite{stoll2017diagonal}. The construction is depicted  in Figure~\ref{fig:construction}.

\begin{lemma}
Let $j \in \{0,1,2,3,4,5\}\setminus \{i\}$ and consider the involution $\tau_j: X_j \mapsto -X_j$ in $\PP^4$. Then 
\[
E_{i,j} = C_i/\langle \tau_j\rangle
\]
is a genus-$1$ curve.
\end{lemma}
\begin{proof}
    The quotient map $\phi: C_i \to E_{i,j}$ has degree $2$. It is ramified at
    $C_i \cap \{X_j =0\} \subset \PP^4$, a set of $8$ points, 
    each with ramification index $2$.
    The Riemann--Hurwitz formula
    gives 
    \(
        2g(C_i) - 2 = 2 \cdot (2g(E_{i,j}) -2) + 8  \cdot (2-1)
    \),
    whence 
    $g(E_{i,j})=1$.
\qed
\end{proof}

We now show how to compute a Weierstrass equation for  $E_{i,j} = C_i/\langle \tau_j\rangle$ by example of the genus-$5$ curve  $C_5 = V(\Omega_0', \Omega_1', \Omega_2') \subset \PP^4$, where
\[
\Omega_0': \sum_{i=0}^{4} X_i^2=0, \quad \Omega_1': \sum_{i=0}^{4} a_i X_i^2=0, \quad \Omega_2': \sum_{i=0}^{4} a_i^2 X_i^2=0.
\]
Moreover, we assume that $j \in \{0,1,2\}$ since the other cases are obtained by permuting the variables. 

First we simplify the equations defining $C_5$ using Gaussian elimination to obtain equations of the form
\[
\begin{matrix}
\Omega_0'':& X_0^2 &&& +\lambda_{0,3}X_3^2& + \lambda_{0,4}X_4^2 &= 0,\\
\Omega_1'':& & X_1^2 &  & +\lambda_{1,3}X_3^2& + \lambda_{1,4}X_4^2 &= 0,\\
\Omega_2'':& & & X_2^2 & +\lambda_{2,3}X_3^2& + \lambda_{2,4}X_4^2 &= 0.
\end{matrix}
\]

The quotient $E_{5,j}$ for $j \in \{0,1,2\}$ is defined as the zero set of the two equations 
\[
\begin{matrix}
\Omega_{j_1}'':& X_{j_1}^2 && +\lambda_{j_1,3}X_3^2& + \lambda_{j_1,4}X_4^2 &= 0,\\
\Omega_{j_2}'':& & X_{j_2}^2   & +\lambda_{j_2,3}X_3^2& + \lambda_{j_2,4}X_4^2 &= 0
\end{matrix}
\]
in $\PP^4$, where $j_1,j_2$ are such that $\{j_1,j_2,j\} = \{0,1,2\}$. This corresponds to the image of $C_5$ under the projection $\pi_j:\PP^5 \to \PP^4$ projecting away from $X_j$. 

Note that $E_{5,j}$ is defined as the intersection of two quadrics in $\PP^3$.
To find a Weierstrass equation for this curve, let $P \in E_{5,j}$ be a rational point. First perform a coordinate transformation such that $P = (0:0:0:1)$ and then consider the projection $\PP^3\to \PP^2$ projecting away from the last coordinate. The restriction of this map to $E_{5,j}$ is birational and in particular the image of $E_{5,j}$ is a curve in $\PP^2$ defined by a cubic equation.

\subsubsection{Singular hyperplane intersections}

\begin{figure}
    \centering
		\begin{tikzpicture}[scale=.5, rotate=90]
		
		
		\draw[thick] (-0.5,0.5) to[quick curve through={(1,1) (2,2) (1,3) (0,4) (1,5) (2,6) (1,7)}]
		(-0.5,7.5);
		\draw[thick] (2.5,0.5) to[quick curve through={(1,1) (0,2) (1,3) (2,4) (1,5) (0,6) (1,7)}]
		(2.5,7.5);
		\clip (-2.5,0) rectangle (4.5,8);
		\end{tikzpicture}
    \caption{Configuration of a singular genus-$5$ curve consisting of two elliptic curves.}
    \label{fig:singulargenus5}
\end{figure}
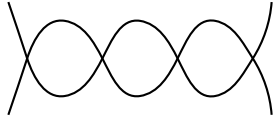

In this part we consider singular curves that arise as the intersection of $\Ksmooth$ with a hyperplane defined as $L:\sum_{i=0}^5 b_i X_i = 0$ for some coefficients $b_i \in k$. Such singular curves have geometric genus $5$ and there are different configurations that can occur.  Since our goal is to find an elliptic curve, we are interested in singular curves that consist of several components with at least one of these an elliptic curve. Here, we discuss the construction of singular curves that consist of two elliptic curves intersecting in $4$ different points. This configuration is depicted in Figure \ref{fig:singulargenus5}.  

Finding parameters $b_i$ such that the intersection is singular can be solved efficiently using linear algebra. For that purpose, one considers the jacobian matrix $M$ of the variety $C = \Ksmooth \cap L$. Let $M(P) \in M_{4,6}(k)$ denote the evaluation of $M$ at a point $P=[x_0: \dots :x_5] \in C$. Then $C$ is singular in $P$ if and only if $\text{rank}(M(P)) = 3$.  Note that the last row of the matrix is given by the vector $b = (b_0, \dots, b_5)$, hence the parameters must be  chosen such that $b$ is a linear combination of the first three rows of the matrix so that $C$ is singular. 

For most choices of $b$, the curve $C$ will consist of only one irreducible component with precisely one singular point. As mentioned before, we intend to construct a curve $C$ with two genus-$1$ components and $4$ singular points. One possibility to achieve this is to choose $b$ such that $b_i=b_j =0$ for two indices $i\neq j$ in $\{0, \dots, 5\}$. In that case, not only $\text{rank}(M(P)) =3$, but $M(P')$ has rank $3$ for every $P' = [x_0': \dots :x_5']$ with $x_k' =x_k$ if $k \notin \{i,j\}$ and $x_k' \in \{\pm x_k\}$ otherwise.

We used this approach for different Kummer surfaces $\Ksmooth$ coming from a superspecial abelian variety. We obtained singular genus-$5$ curves $C$ that consisted of two elliptic curves intersecting in $4$ points. The configuration is depicted in Figure \ref{fig:singulargenus5}.
However, none of the elliptic curves obtained in that way were supersingular.

\subsection{Why do we only obtain ordinary elliptic curves?}
\label{sec:why-ordinary}
In~\S\ref{sec:ungluing_P3} and~\S\ref{sec:ungluing_P5} we succeeded in
constructing elliptic curves from the Kummer surfaces of superspecial
abelian surfaces. However, these elliptic curves were not supersingular in most cases. At first glance this might contradict the intuition that we expect elliptic curves on superspecial abelian surfaces to be supersingular.  To understand this situation, it is necessary to study the preimages of the constructed elliptic curves in the corresponding abelian surface. 

Let us consider the second approach from \S \ref{sec:ungluing_P5}, where we constructed elliptic curves in $\Ksmooth$. If $E \subset \Ksmooth$ is an elliptic curve, then $E' := \phi(E) \subset \Ksing$ is a (possibly singular) genus-$1$ curve.
On the other hand $C = \pi^{-1}(E')$ has genus $1$ if and only if the
cover $\pi$ is unramified along $C$. This means that $E'$ must not go
through any  of the singular points  $S \subset \Ksing$.
The preimages of the sixteen nodes of $\Ksing$ are lines in $\Ksmooth$;
we write $L \subset \Ksmooth$ for this set of lines.
Translating our condition on $E'$ to $\Ksmooth$,
we see that $E$ should not intersect with $L$.
But using explicit descriptions of $L$
(see e.g.~\cite[\S2.2]{castorena2021geometric}),
it is easy to see that there does not exist a hyperplane in $\PP^5$ having trivial intersection with all of these lines. This shows that the elliptic curve $E$ does not correspond to an elliptic curve in $A$.

A similar argument holds for the elliptic curves constructed in \S \ref{sec:ungluing_P3}. The situation in the first approach of \S \ref{sec:ungluing_P5}, where elliptic curves where constructed as quotients of genus-$5$ curves on $\Ksmooth$ is different. One can show that the Jacobian of the genus-$5$ curve $C_i$ as above, is isogenous to  $\prod_{j=1}^5 E_{i,j}$. But it is not clear if there is a relation to $\Jac(Y)$. We leave this as an open question.

\begin{question}
     What is the relation between the elliptic curves $E_{i,j}$ and the Jacobian $\Jac(Y)$ of the initial hyperelliptic curve?
\end{question}

Experimental results show that for each genus-$2$ curve, we find $15$ isomorphism classes of elliptic curves $E_{i,j}$. 
In most cases, the elliptic curves  are not supersingular. 
When starting with $Y: y^2 = x^6-1$, we obtain a mix of ordinary and supersingular elliptic curves.
If it is supersingular, the $j$-invariant is $1728$.


\section{Quantum algorithm for sampling a hard curve}\label{sec:quantum}

On a classical computer, the CGL hash function returns a random curve in the supersingular $\ell$-isogeny graph.  As described in the introduction, if one wishes the curve to be a ``hard curve'' then the drawback to this approach is the need for a trusted party who will throw away the path information generated by the hash function.  Classically, the trusted party seems difficult to avoid.  In this section, we explore the possibility of using a quantum computer to efficiently sample a hard curve from the isogeny graph without leaking any information about the endomorphism ring of the curve.

Although it may be possible to create a quantum algorithm that, when run on a quantum computer, makes the path information inaccessible, there is still a drawback.  Given a curve $E$, we do not know if it was sampled using a classical computer (with an algorithm leaking information about $\End(E)$) or a quantum computer.  Perhaps one can imagine a situation in which all parties inspect the quantum computer and agree it is a quantum computer, and run the program under observation.  However, one may debate whether this situation differs appreciably from the situation in which all parties inspect a classical algorithm designed to delete the path information during its execution, and agree that it will delete it before it can be accessed.  Perhaps one can hope for a means of making the quantum computation ``auditable'' in some way, but we do not have such a method here.  In particular, even if this method samples a uniformly random supersingular curve, it cannot be turned into a hash function in the manner described in the introduction.

Leaving these concerns aside for now, we present below a novel mathematical approach to producing random supersingular curves. We use the idea of continuous-time quantum walks on isogeny graphs of supersingular elliptic curves in characteristic $p$. The idea was first proposed by Kane, Sharif and Silverberg \cite{kane2018quantum,kane2021quantum} for constructing public-key quantum money. In their scheme, quantum walks are carried out over the ideal class group of a quaternion algebra; we adapt these walks to isogeny graphs. The key observation we make here is that the distribution of the curves defined by our sampling algorithm coincides with the limiting distribution of the quantum walks on the graphs.

\subsection{Quantum computing background}

A qubit holds a quantum state that is a superposition (unit length $\C$-linear combination) of the two possible classical states of a bit, i.e. an element of complex norm $1$ of $\C\ket{0} \oplus \C \ket{1}$.  An $n$-qubit quantum register holds a quantum state that is a higher-dimensional analogue: an element $\sum_{x=0}^{2^n-1} \alpha_x \ket{x}$ of complex norm $1$ in $\bigoplus_{0 \le x < 2^n} \C \ket{x}$.  Given any orthonormal basis $\ket{y_i}$ of the $\C$-vector space, we can rewrite the state in that basis: $\sum_i \beta_i \ket{y_i}$. Some of the power of quantum computers comes from the fact that superpositions of $n$ qubits lie in an $2^n$ dimensional state space:  the $n$-fold tensor product of the individual $2$-dimensional state spaces (indeed $(\C\ket{0} \oplus \C \ket{1})^{\otimes n} = \bigoplus_{0 \le x < 2^n} \C \ket{x}$).  Most of those states are \emph{entangled}, meaning that they are not simple tensors in the bases for the individual qubits.

A quantum state $\sum_i \beta_i | y_i \rangle$ cannot be observed except by \emph{measurement in an orthonormal basis $| y_i \rangle$}, a process which collapses the state to one of the basis elements $\ket{y_i}$, where state $\ket{y_i}$ is obtained with probability $\abs{\beta_i}^2$ (the unit length condition implies a valid probability distribution). If there are several registers, we can measure just one, obtaining a superposition of the remaining registers. In a superposition $\sum_{x,y} \alpha_{xy} \ket{x}\ket{y}$, if we measure the first register, we obtain state $C \sum_y \alpha_{x_0y} \ket{x_0}\ket{y}$ (where $C \in \mathbb{R}$ is chosen to scale to unit length) for some $x_0$, with probability $\sum_{y} \abs{\alpha_{x_0 y}}^2$.  

To get started on a quantum computer, one can initialize simple states such as uniform superpositions $\frac{1}{\sqrt{N}} \sum_{i=0}^{N-1} | i \rangle$. A quantum computer then operates on quantum states by unitary operators. Among the most famous is the quantum Fourier transform, whose matrix is that of the inverse discrete Fourier transform.  In particular, it operates by
\[
\sum_{x=0}^{N-1} \alpha_x \ket{x}  \mapsto 
\sum_{x=0}^{N-1} \left( \frac{1}{\sqrt{N}}
\sum_{y=0}^{N-1} \alpha_y e^{2\pi i yx/N}
\right) \ket{x}.
\]
Classical algorithms can be performed in a quantum manner on one quantum register to store the output in another. In particular, for an efficiently computable function $f$ we can perform the operation
\[
\sum_x \alpha_x \ket{x}\ket{0} \mapsto \sum_x \alpha_x \ket{x}\ket{f(x)}.
\]

\subsection{Sampling curves on a quantum computer}

\subsubsection{A na\"ive approach.}
To mimic the CGL algorithm in superposition, we first generate the superposition
\[ \frac{1}{\sqrt{N}} \sum_{x = 1}^N \ket{x}, \]
where $N$ is the number of supersingular curves. Then simultaneously for each $x$, we use the classical CGL algorithm to compute a curve $E_x$, at the end of the path associated to $x$, storing the result in a second register. The resulting superposition is
\[ \frac{1}{\sqrt{N}} \sum_{x = 1}^N \ket{x} \ket{E_x}. \]
Measuring this state collapses the superposition to a classical state $\ket{x}\ket{E_x}$ for some uniformly random $x \in \Z_N$. This is exactly the output of the CGL algorithm for a random input $x$, so the above procedure does not do anything more than the classical CGL.  In particular, the path is stored in the first register.  One way to avoid revealing the path is to apply the quantum Fourier transform to the first register and measure the result. The state we get is
\[ \frac{1}{\sqrt{N}} \sum_{x = 1}^N \omega_N^{xt} \ket{E_x} \]
for some uniformly random $t \in \Z_N$. Now, measuring this state produces a uniformly random curve $E_x$ without revealing anything about the path $x$. However, this approach does not have any advantage over the classical CGL algorithm, as performing the quantum Fourier transform to ``hide" the path information is analogous to including instructions to discard the path information in the classical CGL. In particular, if one measured the first register before the quantum Fourier transform is applied, one could recover the path information. Such runtime interference would not be detectable from the output state alone.

\subsubsection{Continuous-time quantum walk algorithm.}
One way to model random walks on a graph is to apply the adjacency matrix as an operator on the real vector space generated by the vertices (a Markov process).  Na\"ively, one might hope to mimic this on a superposition of the vertices, but unfortunately, this matrix is not unitary. The substitute is the notion of a \emph{quantum walk}, where the adjacency matrix is replaced by its exponential, which is unitary.

The adjacency matrix of the $\ell$-isogeny graph is an $N \times N$ matrix $T_\ell$ called the Brandt matrix. Let us assume, for simplicity, that $T_\ell$ is symmetric.\footnote{This assumption is satisfied for a mild condition on the characteristic $p$.} Let $S$ be the set of supersingular elliptic curves in characteristic $p$. The operator $T_\ell$ acts on the module
\[ M = \bigoplus_{E \in S} \Z E. \]
In the quantum setting, we will work with the complex Euclidean space
\[ \X = M \otimes_\Z \C = \bigoplus_{E \in S} \C \ket{E}.  \]
Note that in order to implement this space on a quantum computer, we use a computational basis of $j$-invariants, so we will include ordinary curves also.  However the random walk, if initiated with a supersingular curve, will restrict itself to the subspace $\X$ generated by the set of supersingular curves.

Let $U_\ell = \exp(i T_\ell)$. The operator $U_\ell$ is unitary (since $T_\ell$ is hermitian) 
and its eigenvalues are $e^{i\lambda}$ for the eigenvalues $\lambda$ of $T_\ell$. The operator $U_\ell^t$ implements a continuous-time quantum walk at time $t$ on the $\ell$-isogeny graph. The application of this for us is that from this quantum walk we can obtain a certain probability distribution on supersingular elliptic curves, and the ability to draw from this distribution to produce a random supersingular elliptic curve (once again, according to this distribution). This is done in the following way: fix an initial supersingular curve $E_0$ and a bound $T > 0$, pick a time $t \in (0, T]$ uniformly at random, compute $U_\ell^t \ket{E_0}$ and measure in the basis $\{ \ket{E} \}_{E \in S}$. The probability of measuring a curve $E \in S$ is then given by \cite[Chapter 16]{QuantumAlgorithms}
\begin{equation}
\label{equ:qwalk-dist}
    p_{E_0 \rightarrow E}(T) = \frac{1}{T} \int_0^T \abs{\braket{E}{e^{iT_\ell t} \vert E_0}}^2 dt.
\end{equation}
For this process to be useful, we must answer two questions about the distribution \eqref{equ:qwalk-dist} on the vertices of the $\ell$-isogeny graph: \begin{enumerate*}[label=(\roman*)] \item How efficient is sampling from this distribution? and \item Do samples leak information about  endomorphism rings?\end{enumerate*}

We comment on the second question first: The question of information leakage requires that we understand the distribution \eqref{equ:qwalk-dist} and the endomorphism rings of its outputs.  However, given an initial curve $E_0$, this distribution seems difficult to analyse.  In particular, it is not the same as the distribution of endpoints of a classical random walk on the $\ell$-isogeny graph.

Regarding efficiency, for any prime $\ell \le \poly(\log N)$, the operator $T_\ell$ is sparse in the sense that there are only $\ell+1 = \poly(\log N)$ nonzero entries in each row or column. Therefore, $T_\ell$ is a good candidate for a Hamiltonian of continuous-time quantum walks;  we can use standard Hamiltonian simulation techniques to implement the quantum walk operator $U_\ell^t$. However, the running time of the best known simulation algorithm depends linearly on $\ell t$ \cite{berry2015hamiltonian}. Therefore, these quantum walks can efficiently be performed only for time $t \le \poly(\log N)$.

\subsubsection{Moving to a limiting distribution.}
To remedy these issues, we consider the limiting distribution of \eqref{equ:qwalk-dist}.  Let $\ket{\phi_j}$, $j = 1, \dots, N$ be a set of eigenvectors of $T_\ell$ and let $\lambda_j$ be the corresponding eigenvalues. It can be shown that \cite[Section 16.6]{QuantumAlgorithms}
\begin{equation}
\label{equ:qwalk-limit}
    \lim_{T \rightarrow \infty} p_{E_0 \rightarrow E}(T) = \sum_{j = 1}^N \abs{\braket{E_0}{\phi_j} \braket{E}{\phi_j}}^2.
\end{equation}
This limiting distribution is more tractable than \eqref{equ:qwalk-dist}, as it is stated in terms of the spectral theory of the graph. 
In practice, for the distribution \eqref{equ:qwalk-dist} to be negligibly close to \eqref{equ:qwalk-limit}, the value $T / (\lambda_j - \lambda_k)$ must be large for any $j, k$. However, the eigenvalues of $T_\ell$ are all in the range $[-2\sqrt{\ell}, 2\sqrt{\ell}]$, so there are some eigenvalues that are exponentially close to each other. This means that for us to assume that we are sampling according to \eqref{equ:qwalk-limit}, we must select $T$ to be exponentially large. But, as mentioned above, we can only implement the walk operator $U_\ell^t$ for polynomially large $t$. Therefore, if we wish to use this nicer distribution, we need a different sampling algorithm which is efficient for larger $T$.

There is a (heuristic) polynomial time algorithm for sampling according to the limiting distribution \eqref{equ:qwalk-limit} using phase estimation. This algorithm is based on the crucial fact that the set of operators $\{ T_\ell \}_{\ell \text{ prime}}$ have a simultaneous set of eigenstates, namely the $\ket{\phi_j}$, $j = 1, \dots, N$ from above. Since $\{ \ket{\phi_j} \}$ is a basis, we can write
\[ \ket{E_0} = \sum_{j = 1}^N \braket{\phi_j}{E_0} \ket{\phi_j}. \]
Now let $\ell_1, \ell_2, \dots, \ell_r$ be a set of primes of size $\poly(\log N)$. Quantum phase estimation is an algorithm to recover the phase (which contains the eigenvalue information) of a unitary operator $U$.  Specifically, if $U \ket{\phi_j} = e^{i\lambda_j} \ket{\phi_j}$ for $j = 1, \ldots, N$, the algorithm recovers an approximation to $\lambda_j$.  We will use phase estimation on the operator $U_{\ell_1}$ with the input state $\ket{E_0}$. Let $\lambda_{1, j}$ be the eigenvalue of $T_{\ell_1}$ corresponding to the eigenstate $\ket{\phi_j}$. Then, because of the relationship between the eigenvalues of $T_\ell$ and those of $U_\ell$, after phase estimation we obtain the state
\begin{equation}
    \label{equ:phase-est}
    \sum_{j = 1}^N \braket{\phi_j}{E_0} \ket{\phi_j} \ket{\tilde{\lambda}_{1, j}}
\end{equation}
where $\abs{\lambda_{1, j} - \tilde{\lambda}_{1, j}} \le 1 / \poly(\log N)$. Measuring the second register (which reveals a value $\tilde{\lambda}_{1,j}$) we obtain a state $\ket{\psi_1}$ that is a projection of the state \eqref{equ:phase-est} onto a smaller subspace $\X_1 \subset \X$. If we repeat this procedure but now with the operator $U_{\ell_2}$ and the input state $\ket{\psi_1}$, we get a new state $\ket{\psi_2}$ that is the projection of $\ket{\psi_1}$ onto a smaller subspace $\X_2 \subset \X_1$. If $r$ is large enough, repeating this procedure for all the remaining $T_{\ell_i}$ we end up with some eigenstate $\ket{\phi_j}$ with probability $\abs{\braket{E_0}{\phi_j}}^2$; see \cite{kane2018quantum,kane2021quantum} for a detailed analysis of this claim. Now, if we measure $\ket{\phi_j}$ in the basis $\{ \ket{E} \}_{E \in S}$, we obtain a curve $E$ with probability $\abs{\braket{E}{\phi_j}}^2$. Therefore, $E$ is a sample from the distribution \eqref{equ:qwalk-limit}.

\subsubsection{Challenges.}
This proposed method still presents a few important questions.  First, a theoretical analysis of the distribution \eqref{equ:qwalk-limit} is needed.  As the $\ell$-isogeny graph is heuristically believed to behave as a random $(\ell+1)$-regular graph, one hopes this distribution will approach the uniform distribution over supersingular curves mod $p$.  Second, the measurement process for phase estimation reveals a series $\tilde{\lambda}_{1,j}, \tilde{\lambda}_{2,j}, \ldots, \tilde{\lambda}_{r, j}$ of approximations to the eigenvalues $\lambda_{1, j}, \lambda_{2, j}, \dots, \lambda_{r, j}$ of the eigenstate $\ket{\phi_j}$ under the operators $T_{\ell_1}, T_{\ell_2}. \dots, T_{\ell_r}$. It is unknown whether revealing this partial eigensystem reveals any information, for example about the likely endomorphism ring of the resulting curve.

\bibliographystyle{cj}
\bibliography{references}

\end{document}